\documentclass[11pt]{amsart}
\usepackage{amssymb,latexsym}
\newdimen\AAdi%
\newbox\AAbo%
\def\AAk#1#2{\s_etbox\AAbo=\hbox{#2}\AAdi=\wd\AAbo\kern#1\AAdi{}}%
\def\AAr#1#2#3{\s_etbox\AAbo=\hbox{#2}\AAdi=\ht\AAbo\raise#1\AAdi\hbox{#3}}%
\font\tenmsb=msbm10 at 11pt \font\sevenmsb=msbm7 at 8pt
\font\fivemsb=msbm5 at 6pt
\newfam\msbfam
\textfont\msbfam=\tenmsb \scriptfont\msbfam=\sevenmsb
\scriptscriptfont\msbfam=\fivemsb

\renewcommand{\theequation}{\thesection\arabic{equation}}

\newcommand{\ba}{\begin{array}}
\newcommand{\ea}{\end{array}}

\newtheorem{theorem}{Theorem}[section]
\newtheorem{lemma}[theorem]{Lemma}
\newtheorem{proposition}[theorem]{Proposition}
\newtheorem{corollary}[theorem]{Corollary}

\theoremstyle{definition}

\theoremstyle{remark}
\newtheorem{remark}[theorem]{Remark}

\begin{document}
\setlength{\baselineskip}{1.2\baselineskip}
\title [Microscopic convexity] {A Microscopic Convexity Principle for
Nonlinear Partial Differential Equations}
\author{Baojun Bian}
\address{Department of mathematics\\
         Tongji University\\Shanghai 200092, People's Republic of
China}
\email{bianbj@mail.tongji.edu.cn}
\author{Pengfei Guan}
\address{Department of Mathematics and Statistics\\
         McGill University\\
         Montreal, Quebec, H3A 2K6, Canada.}
\email{guan@math.mcgill.ca}
\thanks{Research of the first author was
supported in part by NSFC No.10671144 and National Basic Research
Program of China (2007CB814903). Research of the second author was
supported in part by an NSERC Discovery Grant. } \maketitle

\leftskip 0 true cm \rightskip 0 true cm
\renewcommand{\theequation}{\arabic{equation}}
\setcounter{equation}{0} \numberwithin{equation}{section}

\section{Introduction}

Caffarelli-Friedman \cite{CF85} proved a {\it constant rank
theorem} for convex solutions of semilinear elliptic equations in
$\mathbb R^2$, a similar result was
also discovered by Yau \cite{SWYY} at the same time. The result in \cite{CF85}
was generalized to $\mathbb R^n$ by
Korevaar-Lewis \cite{Kl87} shortly after.  This type of
{\it constant rank theorem} is called microscopic convexity
principle. It is a powerful tool in the study of geometric
properties of solutions of nonlinear differential equations, it is
particularly useful in producing convex solutions of differential
equations via homotopic deformations. The great advantage of the
microscopic convexity principle is that it can treat geometric
nonlinear differential equations involving tensors on general
manifolds. The proof of such microscopic convexity principle for
$\sigma_k$-equation on the unit sphere $\mathbb S^n$ by Guan-Ma
\cite{gm} is crucial in the study of the Christoffel-Minkowski
problem. The microscopic convexity principle provides some interesting geometric
properties of solutions to the equation. For symmetric Codazzi tensor, the microscopic
convexity principle yields that the distribution of null space of
the tensor is of constant dimension and it is
parallel. The microscopic convexity principle has been validated
for a varieties of fully nonlinear differential equations
involving the second fundamental forms of hypersurfaces (e.g.,
\cite{gm, GLM1, GMZ, CGM}).

Driven by the pertinent question that under what structural conditions for
partial differential equations so that the microscopic convexity principle
is held, Caffarelli-Guan-Ma \cite{CGM} established such principle for the fully
nonlinear equations of the form:
\begin{eqnarray}
\label{1.1} F(u_{ij}(x)) = \varphi(x, u(x), \nabla u(x)).
\end{eqnarray}
where $F(A)$ is a symmetric and $F(A^{-1})$ is locally convex in
$A$. The similar results were also proved for symmetric tensors on
manifolds in \cite{CGM}, along with several important geometric
applications. It is important to consider equations where $F$
involves other arguments in addition to the Hessian $(u_{ij})$.
For example, it is desirable to include linear elliptic equations
and quasilinear equations with variable coefficients. In many
cases, a solution $v$ to an equation itself may not be convex.
Yet, some of its transformation may be convex (e.g., \cite{BL,
CF85}). If $v$ is a solution of equation (\ref{1.1}), $u=h(v)$ is
a solution of equation
\begin{equation}\label{equ1} F(\nabla^2u,\nabla u,u,x)=0.
\end{equation}
In general, $\nabla^2 u$ may not be separated from the rest of the
arguments. The similar situation also arises in the case of
geometric flow for hypersurfaces.

In this paper, we study the microscopic convexity property for
equation in the form of (\ref{equ1}) and related geometric
nonlinear equations of elliptic and parabolic type. The core for
the microscopic convexity principle is to establish a strong
maximum principle for appropriate constructed functions. The key
is to control ceratin gradient terms of the symmetric tensor to
show that they are vanishing at the end. There have been
significant development of analysis techniques in literature
\cite{CF85, Kl87,gm, GLM1, GMZ, CGM} for this purpose, in
particular the method introduced in \cite{CGM}. They are very
effective to control quadratic terms of the gradient of the
symmetric tensor. For equation (\ref{equ1}), linear terms of such
gradient of symmetric tensor will emerge. All the previous methods
break down for these terms. The main contribution of this paper is
the introduction of new analytic techniques to handle these linear
terms. This type new analysis involves quotients of elementary
symmetric functions near the null set of $\det(u_{ij})$, even
though equation (\ref{equ1}) itself may not be symmetric with
respect to the curvature tensor. The analysis is delicate and has
to be balanced as both symmetric functions in the quotient will
vanish at the null set. This is a novel feature of this paper, it
is another indication that these quotient functions are naturally
embedded with fully nonlinear equations. In a different context,
the importance of quotient functions has been demonstrated in the
beautiful work of Huisken-Sinestrari \cite{HS99}. We believe the
techniques in this paper will find way to solve other problems in
geometric analysis.

To illustrate our main results, we first consider the equations in
flat domain. Let $\Omega$ is a domain in $\mathbb{R}^n $,
$\mathcal{S}^n$ denotes the space of real symmetric $n\times n$
matrices, and $F=F(r,p,u,x)$ is a given function in $
\mathcal{S}^n\times \mathbb{R}^n\times \mathbb{R}\times \Omega $
and elliptic in the sense that
\begin{equation}\label{cond-e}
(\frac{\partial F}{\partial r_{\alpha \beta}}(\nabla^2u,\nabla
u,u,x))>0, \quad \forall x\in \Omega.
\end{equation}

\begin{theorem}\label{microc} Suppose $F=F(r,p,u,x)\in
C^{2,1} (\mathcal{S}^n\times \mathbb{R}^n\times \mathbb{R}\times
\Omega)$ and $F$ satisfies conditions (\ref{cond-e}) and
\begin{equation}\label{cond-c} F(A^{-1},p,u,x) \quad \mbox{is
locally convex in $(A, u, x)$ for each $p$ fixed.}
\end{equation} If $u \in
C^{2,1}(\Omega)$ is a convex solution of (\ref{equ1}), then the
rank of Hessian $(\nabla^2 u(x))$ is constant $l$ in $\Omega$. For
each $x_0\in \Omega$, there exist a neighborhood $\mathcal U$ of
$x_0$ and $(n-l)$ fixed directions $V_1, \cdots, V_{n-l}$ such
that $\nabla^2u(x)V_j=0$ for all $1\le j\le n-l$ and $x\in
\mathcal U$.
\end{theorem}

There is also a parabolic version.

\begin{theorem}\label{microc-para} Suppose $F=F(r,p,u,x,t)\in
C^{2,1} (\mathcal{S}^n\times \mathbb{R}^n\times \mathbb{R}\times
\Omega \times [0,T))$ and $F$ satisfies conditions (\ref{cond-e})
for each $t$ and
\begin{equation}\label{cond-cp} F(A^{-1},p,u,x,t) \quad \mbox{is
locally convex in $(A, u, x)$ for each $(p,t)$ fixed.}
\end{equation}
Suppose $u \in C^{2,1}(\Omega\times [0,T))$ is a convex solution of
the equation
\begin{equation}\label{equ1-parab}
\frac{\partial u}{\partial t}=F(\nabla^2u,\nabla u,u,x, t).
\end{equation}
For each $T>t>0$, let $l(t)$ be the minimal rank of $(\nabla^2
u(x,t))$ in $\Omega$. Then, the rank of $(\nabla^2 u(x,t))$ is
constant for each $T>t>0$ and $l(s)\le l(t)$ for all $s\le t<T$.
For each $0<t\le T$, $x_0\in \Omega$, there exist a neighborhood
$\mathcal U$ of $x_0$ and $(n-l(t))$ fixed directions $V_1,
\cdots, V_{n-l(t)}$  such that $\nabla^2u(x,t)V_j=0$ for all $1\le
j\le n-l(t)$ and $x\in \mathcal U$. Furthermore, for any $t_0\in
[0,T)$, there is $\delta>0$, such that the null space of
$(\nabla^2 u(x,t))$ is parallel in $(x,t)$ for all $x\in \Omega,
t\in (t_0, t_0+\delta)$.
\end{theorem}

An immediate consequence of Theorem \ref{microc} is the validation
of a conjecture raised by Korevaar-Lewis in \cite{Kl87} for convex
solutions of mean curvature type elliptic equation
\begin{equation}\label{4-3}
\sum_{i,j}a^{ij}(\nabla u(x))u_{ij}(x)=f(x,u(x), \nabla u(x))>0.
\end{equation}

\begin{corollary}\label{cor4-1}
Let $\Omega\subset \mathbb R^n$ be a domain. Suppose $u$ is a
convex solution of elliptic equation (\ref{4-3}). If
\begin{equation}\label{qqq-1}
\frac 1{f(x,u,p)} \mbox{is locally convex in $(x,u)$ for each $p$
fixed,}\end{equation} then the Hessian of $u$ is of constant rank
in $\Omega$.
\end{corollary}

Korevaar-Lewis \cite{Kl87} proved that the Hessian of any convex
solution $u$ of elliptic equation (\ref{4-3}) is of constant rank
and $u$ is constant in $n-l$ coordinate directions, provided that
$\frac 1{f(.,p)}$ is strictly convex for any $p$ fixed. They
conjectured that the constant rank result still holds if $\frac
1{f(.,p)}$ is only assumed to be convex. They observed when $n=2$,
this can be deduced from the proofs of Caffarelli-Friedman in
\cite{CF85}. Set
\[F(\nabla^2u, \nabla u, u, x)=-\frac{1}{\sum_{i,j}a^{ij}(\nabla
u(x))u_{ij}(x)}+\frac 1{f(x,u(x), \nabla u(x))},\] Equation
(\ref{4-3}) is equivalent to $F(\nabla^2u, \nabla u, u, x)=0$. It
is straightforward to check that $F$ satisfies Conditions
(\ref{cond-e}) and (\ref{cond-c}) under the assumptions in Corollary
\ref{cor4-1}.

\medskip

We now discuss some geometric equations on general manifolds. Preservation of
convexity is an important issue for the geometric flows of
hypersurfaces (e.g., \cite{Huisken, Andrews} and references
therein). We have the following general result.

\begin{theorem}\label{thmw2-flow}
Suppose $F(A,X, \vec n)$ is elliptic in $A$ and $F(A^{-1},X, \vec
n)$ is locally convex in $(A,X)$ for each fixed $\vec n\in \mathbb
S^n$. Let $M(t)\subset \mathbb R^{n+1}$ be compact hypersurface and it
is a solution of the geometric flow
\begin{equation}\label{flow1}
X_t=-F(g^{-1}h, X, \vec n)\vec n, \ \ t\in (0,T), \ \ M(0)=M_0,
\end{equation}
where $X, \vec n, g, h$ are the position function, outer normal, induced metric
and the second fundamental form of $M(t)$. If $M_0$ is convex,
then $M(t)$ is strictly convex for all $t\in (0,T)$.
\end{theorem}

Alexandrov in \cite{Alex1, Alex2} studied existence and uniqueness of general
nonlinear curvature equations,
\begin{eqnarray}\label{w3.10}
F(g^{-1}h, X, \vec n(X))=0,\quad\forall X \in M,
\end{eqnarray}
where $X$ is the position function of $M$ and $\vec n(X)$ the unit
normal of $M$ at $X$. The following theorem addresses the
convexity problems in \cite{Alex1, Alex2}.

\begin{theorem}\label{thmw2}
Suppose $F(A,X, \vec n)$ is elliptic in $A$ and $F(A^{-1},X, \vec
n)$ is locally convex in $(A,X)$ for each fixed $\vec n\in \mathbb
S^n$. Let $M$ be an oriented immersed connect hypersurface in $
\mathbb R^{n+1}$ with a nonnegative definite second fundamental
form $h$ satisfying equation (\ref{w3.10}), then $h$ is of
constant rank its null space is
parallel. In particular, if $M$ is complete, then there is
$0\le l\le n$ such that $M=M^{l}\times \mathbb R^{n-l}$ for a
strictly convex compact hypersurface $M^l$ in $\mathbb R^{l+1}$.
If in addition $M$ is compact, then $M$ is the boundary of a
strongly convex bounded domain in $\mathbb R^{n+1}$.
\end{theorem}

Theorem \ref{thmw2} shares some similarity with the classical
result of Hartman-Nirenberg in \cite{HN59}.
\medskip

The microscopic convexity principle can be used to prove some
uniqueness theorems in differential geometry in large. An immersed
surface in $\mathbb R^3$ is called Weingarten surface if its
principle curvatures $\kappa_1, \kappa_2$ satisfy relationship
$F(\kappa_1,\kappa_2)=0$ for some function $F$. Alexandrov
\cite{Alex} and Chern \cite{Chern} proved that if $M$ is a closed
convex surface in $\mathbb R^3$ such that $F(\kappa_1,\kappa_2)=0$
for some elliptic $F$ (i.e, $F$ satisfies condition
(\ref{cond-e})), then $M$ is a sphere. In higher dimensions, there
is extensive literature devoted the sphere theorem of immersed
hypersurfaces (e.g., \cite{CY, E-H}). We prove the following
sphere theorem, we refer to \cite{gm, GMZ, CGM} for applications
in classical and conformal geometry, and refer to \cite{GLZ} for
applications in K\"ahler geometry.

\begin{theorem}\label{codazzi-1}
Suppose $(M,g)$ is a compact connected Riemannian manifold of
dimension $n$ with nonnegative sectional curvature, and positive
at one point. Suppose $F(A)$ is elliptic, and $W$ is a Codazzi
tensor on $M$ satisfying equation
\begin{equation}\label{ecod1}
F(g^{-1}W)=0 \quad \text{on $M$.}
\end{equation}
If either
\begin{enumerate} \item $n=2$, or
\item   $n\ge 3$, $W$ is semi-positive definite  and $F(A^{-1})$
is locally convex for $A>0$,\end{enumerate} then $W=cg$ for some
constant $c\ge 0$.
\end{theorem}

Theorem \ref{codazzi-1} was proved by Ecker-Huisken in \cite{E-H}
under the assumption $F$ is concave, we refer Remark \ref{wconds}
for relationship between concavity of $F(A)$ and condition on $F$
in case (2) of Theorem \ref{codazzi-1}. We note that when $n=2$,
only ellipticity assumption on $F$ is needed in Theorem
\ref{codazzi-1}.

\medskip

There is a vast literature devoted to the study of the convexity
of solutions of partial differential equations. There is a theory
of macroscopic nature, where problem is considered in a convex
domain in $\mathbb R^n$ with proper boundary conditions. Korevaar
made breakthroughs in \cite{Ko831, Ko832}, he obtained concavity
maximum principles for a class of quasi-linear elliptic equations
defined convex domains in $\mathbb R^n$ in 1983. His results were
improved by Kennington \cite{Ke85} and by Kawhol \cite{Ka86}. The
theory further developed to its great generality by
Alvarez-Lasry-Lions \cite{All97} in 1997, they established the
existence of convex solution of equation (\ref{equ1}) for state
constraint boundary value under conditions
(\ref{cond-e})-(\ref{cond-c}) and that $F$ satisfies comparison
principle. Microscopic convexity implies macroscopic convexity if
there is a deformation path (e.g., via the method of continuity or
parabolic flow). Theorem \ref{microc} is the microscopic version
of the macroscopic convexity principle in \cite{All97}.

\medskip

The rest of the paper is organized as follows. In section 2, we
introduce a key auxiliary function $q(x)$ and derive certain
negativity properties of this function (Proposition \ref{Prop1-2}
and Corollary \ref{Prop1-1}). In section 3, we establish a strong
maximum principle for function $\phi(x)=\sigma_{l+1}(\nabla^2
u(x))+q(x)$. In section 4, we discuss condition (\ref{cond-c}) and
related results. The last section is devoted to geometric
equations on manifolds.

\medskip

\noindent {\it Acknowledgement:} We would like to thank Professor
Xinan Ma for several helpful discussions. Part of work was done
while the first author was visiting McGill University. He would
like to thank the Department of Mathematics and Statistics at
McGill University for its warm hospitality.

\section{An Auxiliary function}

To establish a microscopic convexity principle, one would like to
prove the rank of $\nabla^2 u$ is of constant rank. It is natural
to consider function $\sigma_{l+1}(\nabla^2 u)$ here $l$ the
minimal rank of $\nabla^2u$. $\nabla^2 u$ is of constant rank is
equivalent to $\sigma_{l+1}(\nabla^2 u)\equiv 0$. It was first
shown by Caffarelli-Friedman in \cite{CF85} that there is a strong
maximum principle for $\sigma_{l+1}(\nabla^2 u)$ when $F=\Delta$
in $\mathbb R^2$. In the subsequential papers \cite{Kl87, gm,
GLM1, GMZ}, this type of maximum principle was establishes for
differential functional $F$ when it is either an elementary
symmetric function of $\nabla^2 u$ or a quotient of them. In these
papers, the analysis relies on the algebraic properties of the
elementary symmetric functions. For general $F$ in (\ref{1.1}),
the test function $\sigma_{l+1}(\nabla^2 u)$ was replaced by
$\sigma_{l+1}(\nabla^2 u)+A\sigma_{l+2}(\nabla^2u)$ ($A$ large).
All these are relied on one special fact: for symmetric function
$F$ in (\ref{1.1}), all the third order derivatives (i.e., the
gradient of the symmetric tensor $\nabla^2u$) which appear in the
process are always in quadratic order. This fact is important for
above mentioned methods to work, we refer Remark \ref{remark3-1}
for a discussion of a unified argument.

When deal with general equation (\ref{equ1}), linear terms of
third order derivatives of $u$ (i.e., the gradient of the
symmetric tensor $\nabla^2u$) will appear. How to control them is
the major challenge. All the test functions considered before
would yield certain "{\it good}" quadratic terms of third order
derivatives which are {\bf not strong} enough for this case, as
linear terms can not be controlled by quadratic terms when they
are assumed to be approaching $0$ (we want prove all of them are
vanishing at the end). We introduce a new auxiliary function which
is composed as a quotient of elementary symmetric functions
$\frac{\sigma_{l+2}(\nabla^2 u)}{\sigma_{l+1}(\nabla^2 u)}$ near
points where $\nabla^2u(x)$ is of minimal rank $l$. Though both
$\sigma_{l+1}(\nabla^2 u)$ and $\sigma_{l+1}(\nabla^2 u)$ vanish
at points where rank of $\nabla^2u(x)$ is $l$, the
Newton-MacLaurine inequality guarantee it is well defined. In
fact, we will show $\frac{\sigma_{l+2}(\nabla^2
u)}{\sigma_{l+1}(\nabla^2 u)}$ has optimal $C^{1,1}$ regularity in
Corollary \ref{Prop1-1}. Furthermore, we will signal out some key
concavity terms of this function in Proposition \ref{Prop1-2} to
dominate the aforementioned linear terms of corresponding third
order derivatives. The quotient function of elementary symmetric
function plays a crucial role in this paper. We also call
attention to the work of \cite{HS99} for some other important
roles of this type of functions in geometric analysis.

\medskip

With the assumptions of $F$ and $u$ in Theorem \ref{microc} and
Theorem \ref{microc-para}, $u$ is automatically in $C^{3,1}$. We will
assume $u \in C^{3,1}(\Omega)$ in the rest of this paper. Let
$W(x)=\nabla^2u(x)$ and $ l=\min_{x\in \Omega} {\rm
rank}(\nabla^2u(x))$. We may assume $l\leq n-1$. Suppose $z_0\in
\Omega$ is a point where $W$ is of minimal rank $l$.

Throughout this paper we assume that $\sigma_{j}(W)=0$ if $j<0$ or
$j>n$. We define for $W=(u_{ij})\in \mathcal S^n$
\begin{eqnarray}\label{def-q}
q(W)= \left \{
  \begin{array}{ll}
   \frac{\sigma_{l+2}(W)}{\sigma_{l+1}(W)},& {\rm if\ \ }
   \sigma_{l+1}(W)>0 \\
   0,& {\rm if\ \ }\sigma_{l+1}(W)=0
   \end{array}
   \right.
\end{eqnarray}
For any symmetric function $f(W)$, we denote
\[
f^{ij}=\frac{\partial f(W)}{\partial u_{ij}},\ \
f^{ij,km}=\frac{\partial^2f(W)}{\partial u_{ij}\partial u_{km}}
\]

For each $z_0\in \Omega$ where $W$ is of minimal rank $l$. We pick
an open neighborhood $\mathcal O$ of $z_0$, for any $x \in
\mathcal O$, let $\lambda_1(x) \le \lambda_2(x) ... \le
\lambda_n(x)$ be the eigenvalues of $W$ at $x$. There is a
positive constant $C>0$ depending only on $\|u\|_{C^{3,1}}$, $W(z_0)$
and $\mathcal O$, such that $\lambda_n(x)\ge \lambda_{n-1}(x) ...
\ge \lambda_{n-l+1}(x) \ge C $ for all $x\in \mathcal O$. Let
$G=\{n-l+1,n-l+2,...,n\}$ and $B=\{1,...,n-l\}$ be the ``good''
and ``bad'' sets of indices respectively. Let
$\Lambda_G=(\lambda_{n-l+1}, ..., \lambda_n)$ be the "good"
eigenvalues of $W$ at $x$ and $\Lambda_B=(\lambda_{1}, ...,
\lambda_{n-l})$ be the "bad" eigenvalues of $W$ at $x$. For the
simplicity, we will also write $G=\Lambda_G$, $B=\Lambda_B$ if
there is no confusion. Note that for any $\delta>0$, we may choose
$\mathcal O$ small enough such that $\lambda_i(x)<\delta$ for all
$i\in B$ and $x\in \mathcal O$.

Set
\begin{equation}\label{def-phi}
\phi=\sigma_{l+1}(W)+q(W)
\end{equation}
where $q$ as in (\ref{def-q}). We will use notation $h=O(f)$ if
$|h(x)|\le C f(x)$ for $x\in \mathcal O$ with positive constant
$C$ under control. It is clear that $\lambda_i=O(\phi)$ for all
$i\in B$.

To get around $\sigma_{l+1}(W)=0$, for $\epsilon>0$ sufficient
small, we consider
\begin{equation}\label{def-qq}
q_{\epsilon}(W)= \frac{\sigma_{l+2}(W_{\epsilon}
)}{\sigma_{l+1}(W_{\epsilon})}, \quad
\phi_{\epsilon}(W)=\sigma_{l+1}(W_{\epsilon})+q_{\epsilon}(W),\end{equation}
where $W_{\epsilon}=W+\epsilon I$. We will also denote
$G_{\epsilon}=(\lambda_{n-l+1}+\epsilon, ...,
\lambda_n+\epsilon)$, $B_{\epsilon}=(\lambda_{1}+\epsilon, ...,
\lambda_{n-1}+\epsilon)$

We will work on $q_{\epsilon}$ to obtain a uniform $C^2$ estimate
independent of $\epsilon$. One may also work directly on $q$ at
the points where $\sigma_{l+1}(\nabla^2 u)\neq 0$ to obtained the
same results in the rest of this section (with all relative
constants independent of chosen point). In any case, we prefer to
work on $q_{\epsilon}$.

Set
\begin{equation}\label{def-v}
v(x)=u(x)+\frac{\epsilon}2|x|^2.\end{equation} We have
$W_{\epsilon}=(\nabla^2 v)$. To simplify the nations, we will
write $q$ for $q_{\epsilon}$, $W$ for $W_{\epsilon}$, $G$ for
$G_{\epsilon}$ and $B$ for $B_{\epsilon}$ with the understanding
that all the estimates will be independent of $\epsilon$. In this
setting, if we pick $\mathcal O$ small enough, there is $C>0$
independent of $\epsilon$ such that
\begin{equation}\label{epi-n1}
\sigma_{l+1}(W(x))\ge C\epsilon, \quad \mbox{ and} \quad  \sigma_1(B(x))\ge C\epsilon,
\quad \mbox{for all $x\in \mathcal O$}. \end{equation}

The importance of the function $q$ is reflected in the following
proposition.

\begin{proposition}\label{Prop1-2}
There are constants $C_1, C_2$ independent of $\epsilon$ such that
at any point $z\in \mathcal O$ with $W$ is diagonal, for any
$\alpha, \beta \in\{1,\cdots,n\}$,
\begin{eqnarray}\label{qqq1}
&& \sum_{i,j,k,m} q^{ij,km}v_{ij\alpha}v_{km\beta} \le
C_1\phi+C_2\sum_{i,j\in B}|\nabla v_{ij}| -2\sum_{i\in B,j\in G}
\frac{\sigma^2_{1}(B|i)-\sigma_{2}(B|i)}{\sigma^2_{1}(B)\lambda_{j}}
v_{ij\alpha}v_{ji\beta} \nonumber \\
&& \quad \quad \quad -\frac{1}{\sigma^3_{1}(B)}\sum_{i\in B}
(\sigma_{1}(B)v_{ii\alpha}-v_{ii}\sum_{j\in
B}v_{jj\alpha})(\sigma_{1}(B)v_{ii\beta}-v_{ii}\sum_{j\in
B}v_{jj\beta})\nonumber
\\ && \quad \quad \quad - \frac{1}{\sigma_{1}(B)}\sum_{i,j\in B, i\neq j}
v_{ij\alpha}v_{ji\beta} -\frac{2}{\sigma^3_{1}(B)}\sum_{i \in B}
v_{ii}\sigma_{1}(B|i)v_{ii\alpha}v_{ii\beta}.
\end{eqnarray}
\end{proposition}
The last three terms in (\ref{qqq1}) will play key role to
dominate linear terms of $v_{ij\alpha}$ ($i,j \in B$) in our proof
of Theorem \ref{microc} in the next section.

\begin{corollary}\label{Prop1-1}\ Let $u \in
C^{3,1}(\Omega)$ be a convex function and $W(x)=(u_{ij}(x)),x\in
\Omega$. Let $ l=\min_{x\in \Omega} {\rm rank}(W(x))$, then the
function $q(x)=q(W(x))$ defined in (\ref{def-q}) is in
$C^{1,1}(\Omega)$.
\end{corollary}

The rest of this section will be devoted to the proof of
Proposition\ref{Prop1-2}, which involves some subtle analysis of
function $q$. The proof of Corollary \ref{Prop1-1} will be given
at the end of this section. In preparation, we will list several
lemmas which are well known. For the sack of completeness, we will
provide the proofs. Suppose $W$ is any $n\times n$ diagonal
matrix, we denote $(W|i)$ to be the $(n-1)\times (n-1)$ matrix
with $i$th row and $i$th column deleted, and denote $(W|ij)$ to be
the $(n-2)\times (n-2)$ matrix with $i,j$th rows and $i,j$th
columns deleted.

\begin{lemma}\label{Lemma1-1}\ Suppose $W$ is diagonal.
Then we have
\begin{eqnarray*}
q^{ij}=\left \{
\begin{array}{ll}
      \frac{\sigma_{l+1}(W)\sigma_{l+1}(W|i)
      -\sigma_{l+2}(W)\sigma_{l}(W|i)}{\sigma^2_{l+1}(W)},
          &   {\rm if\ }i=j   \\
     0    &   {\rm if\ }i\neq j
     \end{array}, and
   \right.
\end{eqnarray*}
(a). if $i=m,j=k,i\neq j$, then
\begin{eqnarray*}
q^{ij,km} = \frac{\sigma_{l}(W|ij)}{\sigma_{l+1}(W)}+
    \frac{\sigma_{l+2}(W)\sigma_{l-1}(W|ij)}{\sigma^2_{l+1}(W)}
\end{eqnarray*}
(b). if $i=j=k=m$, then
\begin{eqnarray*}
q^{ij,km}=
-2\frac{\sigma_{l}(W|i)}{\sigma^3_{l+1}(W)}[\sigma_{l+1}(W)\sigma_{l+1}(W|i)
    -\sigma_{l}(W|i)\sigma_{l+2}(W|i)]
\end{eqnarray*}
(c). if $i=j,k=m,i\neq k$, then
\begin{eqnarray*}
q^{ij,km}= \frac{\sigma_{l}(W|ik)}{\sigma_{l+1}(W)}-
\frac{\sigma_{l+1}(W|i)\sigma_{l}(W|k)}{\sigma^2_{l+1}(W)}
     - \frac{\sigma_{l+1}(W|k)\sigma_{l}(W|i)}{\sigma^2_{l+1}(W)}
\end{eqnarray*}
\begin{eqnarray*}
-\frac{\sigma_{l+2}(W)\sigma_{l-1}(W|ik)}{\sigma^2_{l+1}(W)} +2
\frac{\sigma_{l+2}(W)\sigma_{l}(W|i)\sigma_{l}(W|k)}{\sigma^3_{l+1}(W)}
\end{eqnarray*}
(d). otherwise
\begin{eqnarray*}
q^{ij,km}=     0
\end{eqnarray*}
\end{lemma}

\noindent \textbf{Proof.}\ Since $W$ is diagonal, it follows from
Proposition 2.2 in \cite{gm}
\begin{eqnarray*}
\frac{\partial \sigma_{\gamma}(W)}{\partial v_{ij}}= \left \{
  \begin{array}{ll}
    \sigma_{\gamma -1}(W|i),    &   {\rm if\ }i=j  \nonumber \\
    0,    &   {\rm if\ }i\neq j
    \end{array}
   \right.
\end{eqnarray*}
and
\begin{eqnarray*}
\frac{\partial^2 \sigma_{\gamma}(W)}{\partial v_{ij} \partial
v_{km}}= \left \{
  \begin{array}{ll}
     \sigma_{\gamma -2}(W|ik),    &   {\rm if\ }i=j, k=m, i\neq k \nonumber  \\
     -\sigma_{\gamma -2}(W|ij),    &   {\rm if\ }i=m,j=k,i\neq j \nonumber \\
     0,    &   {\rm otherwise}
     \end{array}
   \right.
\end{eqnarray*}
for $1\leq \gamma \leq n$. We obtain thus
\[
    \sigma_{l+1}^{ij}=\frac{\partial \sigma_{l+1}}{\partial W_{ij}}= \left \{
  \begin{array}{ll}
    \sigma_{l}(W|i),    &   {\rm if\ }i=j   \\
     0,    &   {\rm if\ }i\neq j
     \end{array}
   \right.
\]
and
\begin{equation}\label{1-3}
    \sigma_{l+1}^{ij,km}=\frac{\partial^2 \sigma_{l+1}}{\partial W_{ij} \partial W_{km}}
    = \left \{
  \begin{array}{ll}
     \sigma_{l-1}(W|ik),    &   {\rm if\ }i=j, k=m, i\neq k   \\
     -\sigma_{l-1}(W|ij)    &   {\rm if\ }i=m,j=k,i\neq j \\
     0    &   {\rm otherwise}
     \end{array}
   \right.
\end{equation}

A direct computation yields
\begin{equation}\label{1-4}
q^{ij} =\frac{1}{\sigma_{l+1}(W)}\frac{\partial
\sigma_{l+2}(W)}{\partial
v_{ij}}-\frac{\sigma_{l+2}(W)}{\sigma^2_{l+1}(W)}\frac{\partial
\sigma_{l+1}(W)}{\partial v_{ij}}
\end{equation}
and
\[
 q^{ij,km}=\frac{1}{\sigma_{l+1}(W)}\frac{\partial^2
\sigma_{l+2}(W)}{\partial v_{ij}\partial v_{km}}
-\frac{1}{\sigma^2_{l+1}(W)} \frac{\partial
\sigma_{l+2}(W)}{\partial v_{ij}} \frac{\partial
\sigma_{l+1}(W)}{\partial v_{km}}
\]
\[
 -\frac{1}{\sigma^2_{l+1}(W)}
\frac{\partial \sigma_{l+2}(W)}{\partial v_{km}} \frac{\partial
\sigma_{l+1}(W)}{\partial v_{ij}}
-\frac{\sigma_{l+2}(W)}{\sigma^2_{l+1}(W)}\frac{\partial^2
\sigma_{l+1}(W)}{\partial v_{ij}\partial v_{km}}
\]
\begin{equation}\label{1-5}
+2\frac{\sigma_{l+2}(W)}{\sigma^3_{l+1}(W)}\frac{\partial
\sigma_{l+1}(W)}{\partial v_{ij}}\frac{\partial
\sigma_{l+1}(W)}{\partial v_{km}}
\end{equation}

The lemma follows from (\ref{1-4}) and (\ref{1-5}). \qed

\medskip

\begin{lemma}\label{Lemma1-2}\ Suppose $W$ is diagonal, then
\[
q^{ij} =\left \{
\begin{array}{ll}
\frac{\sigma^2_{1}(B|i)-\sigma_{2}(B|i)}{\sigma^2_{1}(B)}+O(\phi),
&   {\rm if\ }i=j\in B  \\
O(\phi),    &   {\rm if\ }i=j\in G  \\
 0,& {\rm if\ } i\neq j.
     \end{array}
   \right.
\]
Furthermore $q^{ij,km}$ can be computed as follows:
\begin{enumerate}
\item If $i,j,k,m\in G$,  \begin{eqnarray*} q^{ij,km}=O(\phi)
\end{eqnarray*}
\item If $j\in G, i\in B$, \begin{eqnarray*}
q^{ji,ij}=q^{ij,ji}=-\frac{\sigma^2_{1}(B|i)-\sigma_{2}(B|i)}{\sigma^2_{1}(B)v_{jj}}+O(\phi)
\end{eqnarray*}
\item If $i,j\in B, i\neq j$,
\begin{eqnarray*}
q^{ij,ji}=-\frac{1}{\sigma_{1}(B)}+O(1)
\end{eqnarray*}
\item If $i \in B$,
\begin{eqnarray*}
q^{ii,ii}=-\frac{2}{\sigma^3_{1}(B)}(\sigma_{1}(B)\sigma_{1}(B|i)-\sigma_{2}(B|i))+O(1)
\end{eqnarray*}
\item If $i\in B, k\in G$,
\begin{eqnarray*}
q^{kk,ii}=q^{ii,kk}=O(1)
\end{eqnarray*}
\item If $i,k\in B, i\neq k$,
\begin{eqnarray*}
q^{ii,kk}=\frac{2\sigma_{2}(B)-\sigma^2_{1}(B)+(v_{ii}+v_{kk})
\sigma_{1}(B)}{\sigma^3_{1}(B)}+O(1)
\end{eqnarray*}
\item otherwise
\begin{eqnarray*}
q^{ij,km}=     0.
\end{eqnarray*}
\end{enumerate}
\end{lemma}

\noindent \textbf{Proof.}\ From \cite{gm} we conclude that for
$W=(G,B)$ and $\gamma \geq l$,
\[
\sigma_{\gamma}(W)=\sum_{k=0}^{l}
\sigma_{k}(G)\sigma_{\gamma-k}(B),
\]
and
\[
\sigma_{\gamma}(W|i)=\sum_{k=0}^{l}
\sigma_{k}(G)\sigma_{\gamma-k}(B|i), \quad \mbox{for $i\in B$};
\]

\[
\sigma_{\gamma}(W|i)=\sum_{k=0}^{l-1}
\sigma_{k}(G|i)\sigma_{\gamma-k}(B), \quad \mbox{ for $i\in G$}:\]
\[
\sigma_{\gamma}(W|ij)=\sum_{k=0}^{l-2}
\sigma_{k}(G|ij)\sigma_{\gamma-k}(B), \quad \mbox{ for $i,j\in
G$};\]
\[
\sigma_{\gamma}(W|ij)=\sum_{k=0}^{l-1}
\sigma_{k}(G|i)\sigma_{\gamma-k}(B|j), \quad \mbox{ for $i\in
G,j\in B$}\]
\[
\sigma_{\gamma}(W|ij)=\sum_{k=0}^{l}
\sigma_{k}(G)\sigma_{\gamma-k}(B|ij), \quad \mbox{ for $i,j\in
B$},\] where $\sigma_{\gamma-k}(B)=0$ if $\gamma-k > n-l$. The
lemma follows directly from lemma \ref{Lemma1-1} and above
formulae. \qed

\medskip

Next we establish an estimate for third order derivatives of
convex functions.

\begin{lemma}\label{Lemma1-3}  Assume $u \in
C^{3,1}(\Omega)$ is a convex function. Then there exists a
positive constant C depending only on $dist\{\mathcal O, \partial
\Omega\}$ and $\|v\|_{C^{3,1}(\Omega)}$ such that
\begin{equation}\label{eq1-7}
|v_{ij \alpha}(x)|\leq C\Big(\sqrt{v_{ii}(x)}+\sqrt{v_{jj}(x)} \Big)
\end{equation}
for all $x\in \mathcal O$ and $1\leq i,j,\alpha\leq n$.\end{lemma}

\noindent \textbf{Proof.}\  It follows from convexity of $v$ that
for any direction $\eta \in R^n$ with $|\eta|=1$
\[
v_{\eta\eta}(x)\geq 0
\]
for all $x\in \Omega$. It's well known that for any nonnegative
$C^{1,1}$ function $h$, $|\nabla h(x)|\le C h^{\frac 12}(x)$ for all
$x\in \mathcal O$, where $C$ depending only on
$\|h\|_{C^{1,1}(\Omega)}$ and $dist\{\mathcal O,
\partial \Omega\}$ (e.g., see \cite{T}). We now infer
\[
|v_{\eta \eta \alpha}(x)|\leq C\sqrt{v_{\eta \eta}(x)}.
\]
where $C$ is a positive constant depending only on $dist\{\mathcal
O,
\partial \Omega\}$ and $\|v_{\eta\eta}\|_{C^{1,1}(\Omega)}$
(which can be controlled by $\|u\|_{C^{3,1}(\Omega)}$).  Now set
$\eta = i$ if $i=j$ and
\[
  \eta  =
    \frac{1}{\sqrt{2}}(e_i+e_j) \quad   {\rm if} \quad i\neq j.
\]
Proof of Lemma \ref{Lemma1-3} is complete. \qed

\medskip

\begin{remark}\label{remark3-1} In \cite{CGM}, test function
$\phi(x)=\sigma_{l+1}(\nabla^2 u(x))+A \sigma_{l+2}(\nabla^2
u(x))$ was introduced. The term $A \sigma_{l+2}(\nabla^2 u(x))$
was used there to overcome quadratic terms of the third order
derivatives. With Lemma \ref{Lemma1-3}, these quadratic terms of
the third order derivatives in fact can be controlled by
$\sigma_{l+1}(\nabla^2 u(x))$. Therefore, all the arguments in
\cite{CGM} can carry through for simpler test function
$\phi(x)=\sigma_{l+1}(\nabla^2 u(x))$. Nevertheless, for general
equation (\ref{equ1}), we will see in the next section that linear
terms of the third order derivatives will appear, the auxiliary
function $q(x)$ will play crucial role to control these terms.
\end{remark}

\bigskip

\noindent \textbf{Proof of Proposition \ref{Prop1-2}.}\ Let us
divide ${\sum}_{i,j,k,m} q^{ij,km}v_{ij\alpha}v_{km\beta}$ into
three parts according to Lemma \ref{Lemma1-1}:
\begin{equation}\label{q-div}
\sum_{i,j,k,m}q^{ij,km}(W(z))v_{ij\alpha}v_{km\beta}
=I_{\alpha\beta}+II_{\alpha\beta}+III_{\alpha\beta},
\end{equation}
where
\[
I_{\alpha\beta}=\sum_{i\neq j}q^{ij,ji}v_{ij\alpha}v_{ji\beta},
\]
\[
II_{\alpha\beta}=\sum_{i=1}^n q^{ii,ii}v_{ii\alpha}v_{ii\beta}
\]
and
\[
III_{\alpha\beta}=\sum_{i\neq k}q^{ii,kk} v_{ii\alpha}v_{kk\beta}.
\]
Lemma \ref{Lemma1-2} yields
\begin{eqnarray}\label{1-9}
I_{\alpha\beta}&=&(\sum_{i,j\in G,i\neq j}+\sum_{i\in B, j\in G}
+\sum_{j\in B, i\in G}+\sum_{i,j\in B,i\neq j})
q^{ij,ji}v_{ij\alpha}v_{ji\beta} \nonumber \\
& &= O(\phi)+O(\sum_{i,j\in B}|\nabla v_{ij}|)-
\frac{1}{\sigma_{1}(B)}\sum_{i,j\in B, i\neq j}
v_{ij\alpha}v_{ji\beta} \nonumber \\
& &-2\sum_{i\in B, j\in G}\frac{\sigma_1^2(B|i)
-\sigma_{2}(B|i)}{\sigma^2_{1}(B)v_{jj}}v_{ij\alpha}v_{ji\beta}.
\end{eqnarray}
It follows that from Lemma \ref{Lemma1-2}
\begin{eqnarray}\label{1-11}
II_{\alpha\beta}&=&(\sum_{i\in G} + \sum_{i\in B})
q^{ii,ii}v_{ii\alpha}v_{ii\beta} \nonumber
\\
&= & O(\phi)+O(\sum_{i,j\in B}|\nabla v_{ij}|)-2\sum_{i \in B}
\frac{\sigma_{1}(B)\sigma_{1}(B|i)-\sigma_{2}(B|i)}{\sigma^3_{1}(B)}v_{ii\alpha}v_{ii\beta}
\end{eqnarray}
and
\begin{eqnarray}\label{1-12}
III_{\alpha\beta}&=&(\sum_{i,j\in G,i\neq j}+\sum_{i\in B, j\in G}
+\sum_{j\in B, i\in G}+\sum_{i,j\in B,i\neq j})
q^{ii,jj}v_{ii\alpha}v_{jj\beta} \nonumber \\
&=& O(\phi)+O(\sum_{i,j\in B}|\nabla v_{ij}|)+ \sum_{i\neq
j,i,j\in B }
\frac{2\sigma_{2}(B)-\sigma^2_{1}(B)+(v_{ii}+v_{jj})\sigma_{1}(B)}{\sigma^3_{1}(B)}
v_{ii\alpha}v_{jj\beta}.
\end{eqnarray}
By the identity, for any indices set $A$,
\begin{eqnarray}\label{id1}
&& \sum_{i,j\in A, i\neq j}
[2\sigma_{2}(A)-\sigma^2_{1}(A)+(v_{ii}+v_{jj})\sigma_{1}(A)]
v_{ii\alpha}v_{jj\beta}\nonumber \\
&& \quad \quad \quad \quad -2\sum_{i \in A}
[\sigma_{1}(A)\sigma_{1}(A|i)-\sigma_{2}(A|i)]v_{ii\alpha}v_{ii\beta}
\nonumber
\\
&&  \quad \quad =-\sum_{i\in
A}(\sigma_{1}(A)v_{ii\alpha}-v_{ii}\sum_{j\in A}v_{jj\alpha})
(\sigma_{1}(A)v_{ii\beta}-v_{ii}\sum_{j\in A}v_{jj\beta})\nonumber
\\
&& \quad \quad \quad \quad - 2\sum_{i \in A}
v_{ii}\sigma_{1}(A|i)v_{ii\alpha}v_{ii\beta}.
\end{eqnarray}
In particular, setting $A=B$ in (\ref{id1}), we deduce
\begin{eqnarray}\label{1-13}
II_{\alpha\beta}+III_{\alpha\beta}&=& O(\phi)+O(\sum_{i,j\in
B}|\nabla v_{ij}|)-\frac{2}{\sigma^3_{1}(B)}\sum_{i \in B}
v_{ii}\sigma_{1}(B|i)v_{ii\alpha}v_{ii\beta} \nonumber \\
&-& \frac{1}{\sigma^3_{1}(B)}\sum_{i\in
B}(\sigma_{1}(B)v_{ii\alpha}-v_{ii}\sum_{j\in B}v_{jj\alpha})
(\sigma_{1}(B)v_{ii\beta}-v_{ii}\sum_{j\in B}v_{jj\beta}).
\end{eqnarray}
\qed

\medskip

Finally, we prove Corollary \ref{Prop1-1}.

\noindent \textbf{Proof of Corollary \ref{Prop1-1}.}\ We only need
to consider a small neighborhood $\mathcal O$ of these point $p\in
\Omega$ such that the minimal rank is attained at $p$. For such
fixed point $z\in \mathcal O$, we may assume $W(z)$ is diagonal by
a rotation. We thus obtain for any fixed $\alpha $ and $\beta$
\begin{equation}\label{eq1-8}
\frac{\partial^2 q(z)}{\partial x_{\alpha}\partial
x_{\beta}}=\sum_{i,j}q^{ij}(W(z))u_{ij\alpha\beta}
+\sum_{i,j,k,m}q^{ij,km}(W(z))u_{ij\alpha}u_{km\beta}
\end{equation}

Since
$0\leq\frac{\sigma^2_{1}(B|i)-\sigma_{2}(B|i)}{\sigma^2_{1}(B)}\leq
1$, by Lemma \ref{Lemma1-2}
\[
 |q^{ij}(W(z))|\leq C
 \]
for some constant $C$ under control. It yields the estimate for
the first term in (\ref{eq1-8})
\[
\|q^{ij}(W(z))u_{ij\alpha\beta}\|\leq C\|u\|_{C^{3,1}(\Omega)}\leq
C
\]

We treat the second term in (\ref{eq1-8}). By Lemma
\ref{Lemma1-3}, for $i,j\in B$
\begin{equation}\label{1-10}
|u_{ij\alpha}|\leq C(\sqrt{u_{ii}(x)}+\sqrt{u_{jj}(x)})\leq
C\sqrt{\sigma_{1}(B)}.
\end{equation}
Noting that $u_{jj}\geq C>0,j\in G$ and
$0\leq\frac{\sigma^2_{1}(B|i)-\sigma_{2}(B|i)}{\sigma^2_{1}(B)}\leq
1$. It now follows from Proposition \ref{Prop1-2},
\[
|\frac{\partial^2 q(W(z))}{\partial x_{\alpha}\partial
x_{\beta}}|\leq C
\]
for all $z\in \mathcal O$.\qed

\section{A strong maximum principle}

In this section, we prove a strong maximum principle for $\phi$
defined in (\ref{def-phi}) for equation (\ref{equ1}). We may prove
the same result for equation (\ref{equ1-parab}) and make Theorem
\ref{microc} as a corollary of Theorem \ref{microc-para}. But we
prefer to work on elliptic case first. The parabolic version will
be proved at the end of next section with some minor modification.

We denote $\mathcal{S}^{n}$ to be the set of all real symmetric
$n\times n$ matrices, and denote $\mathcal{S}^{n}_{+}\subset
\mathcal{S}^{n}$ to be the set of all positive definite symmetric
$n\times n$ matrices. Let $\mathbb O_{n}$ be the space consisting
all $n\times n$ orthogonal matrices. We define
\begin{eqnarray*} \mathcal{S}_{n-1}=\{
Q\left (
\begin{array}{cc}
       0 & 0 \\ 0 & B  \end{array} \right )Q^T \quad |\quad
       \mbox{$\forall Q\in \mathbb O_n$, $\forall B\in \mathcal{S}^{n-1}$
       }\},
       \end{eqnarray*}
and for given $Q\in \mathbb O_{n}$,
\begin{eqnarray*} \mathcal{S}_{n-1}(Q)=\{
Q\left (
\begin{array}{cc}
       0 & 0 \\ 0 & B  \end{array} \right )Q^T \quad |\quad
       \mbox{ $\forall B\in \mathcal{S}^{n-1}$ }\}.
       \end{eqnarray*}
Therefore $\mathcal{S}_{n-1}, \mathcal{S}_{n-1}(Q) \subset
\mathcal{S}^{n}$. For any function $F(r,p,u,x)$, we denote
\begin{eqnarray}\label{eq2-1}
& &F^{\alpha\beta}=\displaystyle\frac{\partial F}
 {\partial r_{\alpha\beta}},\ \ F^{u}=\displaystyle\frac{\partial F}
 {\partial u},\ \ F^{x_i}=\displaystyle\frac{\partial F}
 {\partial x_i},\ \ F^{\alpha\beta,\gamma\eta}=\displaystyle\frac{\partial^2 F}
 {\partial r_{\alpha\beta}\partial r_{\gamma\eta}}, \ \ F^{\alpha \beta,u}
 =\displaystyle\frac{\partial^2
 F}
 {\partial r_{\alpha\beta}\partial u}, \nonumber \\
& & \quad F^{\alpha \beta,x_k}=\displaystyle\frac{\partial^2
F}{\partial r_{\alpha\beta}\partial x_k}, \ \
F^{u,u}=\displaystyle\frac{\partial^2 F}
 {\partial^2 u}, \ \ F^{u,x_i}=\displaystyle\frac{\partial^2 F}
 {\partial u\partial x_i}, \ \ F^{x_i,x_j}=\displaystyle\frac{\partial^2 F}
 {\partial x_i\partial x_j}.
\end{eqnarray}
For any $p$ fixed and $Q\in \mathbb O_{n}$, $(A,u,x)\in \mathcal
S_{n-1}(Q) \times \mathbb R\times \mathbb R^n$, we set
\[ X^*_F=((F^{\alpha
\beta}(A,p,u,x)), -F^u(A,p,u,x), -F^{x_1}(A,p,u,x), \cdots,
-F^{x_1}(A,p,u,x))\] as a vector in $\mathcal S^n\times \mathbb
R\times \mathbb R^n$. Set
\begin{eqnarray} \Gamma^{\bot}_{X^*_F}=\{\tilde X\in \mathcal{S}_{n-1}(Q)\times \mathbb R\times
\mathbb R^n \quad | \quad <\tilde X, X^*_F>=0\},\end{eqnarray}

Let $B\in \mathcal{S}^{n-1}_+, A=B^{-1}$ and
\[
\tilde{B}=\left ( \begin{array}{cc}
       0 & 0 \\ 0 & B  \end{array} \right ),\ \ \tilde{A}=\left ( \begin{array}{cc}
       0 & 0 \\ 0 & A  \end{array} \right ).
\]
For any given $Q \in \mathbb O_{n}$ and $\tilde X=((X_{ij}), Y,
Z_1,\cdots,Z_n) \in \mathcal{S}_{n-1}(Q)\times \mathbb R\times
\mathbb R^n$, we define a quadratic form
\begin{eqnarray}\label{eq2-5}
Q^*(\tilde X, \tilde X)&=&\sum_{i,j,k,l=1}^n F^{ij,kl}X_{ij}X_{kl}
+2\sum_{i,j,k,l=1}^n F^{ij}(Q\tilde{A}Q^T)_{kl}
X_{ik}X_{jl}+\sum_{i,j=1}^n
F^{x_i,x_j}Z_iZ_j \nonumber \\
&&-2\sum_{i,j=1}^n F^{ij,u}X_{ij}Y -2\sum_{i,j,k=1}^n
F^{ij,x_k}X_{ij}Z_k+ 2\sum_{i=1}^n F^{u,x_i}YZ_i+ F^{u,u}Y^2,
\end{eqnarray}
where functions $F^{ij,kl}, F^{ij}, F^{u,u}, F^{ij,u}, F^{ij,x_k},
F^{u,x_i}, F^{x_i,x_j}$ are evaluated at  $(Q\tilde{B}Q^T,p,u,x)$.

We first state a lemma, it's proof will be given in next section
(after Corollary \ref{corollary2-2}).
\begin{lemma}\label{lemma-cdc}
If $F$ satisfies condition (\ref{cond-c}), then for each $p\in
\mathbb R^n$,
\begin{eqnarray}\label{wwcond-c}
F(0,p,u,x) \quad \mbox{is locally convex in $(u,x)$, and }
Q^*(\tilde X, \tilde X) \ge 0, \forall \tilde X \in \Gamma^{\bot}_{
X^*_F}.
\end{eqnarray}\end{lemma}

The following theorem is the core of this paper. Theorem
\ref{microc} is a direct consequence of Theorem \ref{microc-w} and
Lemma \ref{lemma-cdc}.

\begin{theorem}\label{microc-w} Suppose that the function $F$ satisfies
conditions (\ref{cond-e}) and (\ref{wwcond-c}), let $u \in
C^{3,1}(\Omega)$ is a convex solution of (\ref{equ1}). If
$\nabla^2 u$ attains minimum rank $l$ at certain point $x_0\in
\Omega$, then there exist a neighborhood $\mathcal O$ of $x_0$ and
a positive constant $C$ independent of $\phi$ (defined in
(\ref{def-phi})), such that
\begin{equation}\label{2-6}
 \sum_{\alpha, \beta} F^{\alpha\beta}\phi_{\alpha\beta}(x)\leq
 C(\phi(x) +|\nabla\phi(x)|), \quad \forall x\in \mathcal O.
\end{equation}
In turn, $\nabla^2 u$ is of constant rank in $\mathcal O$. Moreover,
for each $x_0\in \Omega$, there exist a neighborhood $\mathcal U$ of
$x_0$ and $(n-l)$ fixed directions $V_1, \cdots, V_{n-l}$ such
that $\nabla^2u(x)V_j=0$ for all $1\le j\le n-l$ and $x\in
\mathcal U$.
\end{theorem}

\noindent \textbf{Proof of Theorem \ref{microc-w}.}\
Let $u\in C^{3,1}(\Omega)$ be a convex solution of equation
(\ref{equ1}) and $W(x)=(u_{ij}(x))$. For each $z_0\in \Omega$
where $W=(\nabla^2u)$ attains minimal rank $l$. We may assume
$l\le n-1$, otherwise there is nothing to prove. As in the
previous section, we pick an open neighborhood $\mathcal O$ of
$z_0$, for any $x \in \mathcal O$, let $G=\{n-l+1,n-l+2,...,n\}$
and $B=\{1,...,n-l\}$ be the ``good'' and ``bad'' sets of indices
for eigenvalues of $\nabla^2u(x)$ respectively.

Setting $\phi$ as (\ref{def-phi}), then we see from Corollary
\ref{Prop1-1} that $\phi\in C^{1,1}(\mathcal O)$ ,
\[
\phi(x)\geq 0,\ \phi(z_0) = 0
\]
and there is a constant $C>0$ such that for all $x\in \mathcal O$,
\[
\frac 1C \sigma_{1}(B)(x)\le \phi(x) \le C\sigma_{1}(B)(x), \
\frac 1C \sigma_{1}(B)(x)\le \sigma_{l+1}(x)\le C\sigma_{1}(B)(x).
\]

We shall fix a point $z\in \mathcal O$
and prove (\ref{2-6}) at $z$. For each $z\in \mathcal O$ fixed,
letting $\lambda_1 \le \lambda_2 ... \le \lambda_n$ be the
eigenvalues of $W(z)=(u_{ij}(z))$ at $z$, we can rotate coordinate
so that $W(z)=(u_{ij}(z))$ is diagonal, and $u_{ii}(z)=\lambda_i,
i=1,\cdots, n$. We note that all quantities involving $g,q$ and
$\phi$ are invariant under rotation.

Again, as in the previous section, we will avoid to deal with
$\sigma_{l+1}(W)=0$ by considering for $W_{\epsilon}$ (defined in
(\ref{def-qq})) for $\epsilon>0$ sufficient small, with
$W_{\epsilon}=W+\epsilon I$,
$G_{\epsilon}=(\lambda_{n-l+1}+\epsilon, ...,
\lambda_n+\epsilon)$, $B_{\epsilon}=(\lambda_{1}+\epsilon, ...,
\lambda_{n-1}+\epsilon)$. We note that $W_{\epsilon}$ is the
Hessian of function $u_{\epsilon}(x)=u(x)+\frac{\epsilon}2 |x|^2$.
This function $u_{\epsilon}(x)$ satisfies equation
\begin{equation}\label{eq-v1}
F(\nabla^2 u_{\epsilon}, \nabla  u_{\epsilon}, u_{\epsilon},
x)=R_{\epsilon},
\end{equation}
where $R_{\epsilon}(x)=F(\nabla^2 u_{\epsilon}, \nabla
u_{\epsilon}, u_{\epsilon}, x)-F(\nabla^2u, \nabla u, u, x)$.
Since $u\in C^{3,1}$, we have
\begin{equation}\label{eq-RRR}
|R_{\epsilon}(x)|\le C\epsilon, \quad |\nabla R_{\epsilon}(x)|\le
C\epsilon,\quad |\nabla^2R_{\epsilon}(x)|\le C\epsilon, \quad
\forall x\in \mathcal O.\end{equation}

We will work on equation (\ref{eq-v1}) to obtain differential
inequality (\ref{2-6}) for $\phi_{\epsilon}$ defined in
(\ref{def-qq}) with constant $C_1, C_2$ independent of $\epsilon$.
Theorem \ref{microc-w} would follow by letting $\epsilon \to 0$.

Set $v=u_{\epsilon}$, in the rest of this section, we will write
$W$ for $W_{\epsilon}$, $G$ for $G_{\epsilon}$,  $B$ for
$B_{\epsilon}$, $q$ for $q_{\epsilon}$ and $\phi$ for
$\phi_{\epsilon}$, with the understanding that all the estimates
will be independent of $\epsilon$. We note that by (\ref{epi-n1}),
we have
\begin{equation}\label{epi-n2}
 \epsilon\le C\phi(x), \quad \mbox{for all $x\in \mathcal O$,}
\end{equation}
and $v$ satisfies equation
\begin{equation}\label{eq-v} F(\nabla^2 v, \nabla v, v, x)=R(x), \end{equation}
with $R(x)$ under control as follows,
\begin{equation}\label{R-est1}
|\nabla^j R(x)|\le C\phi(x), \quad \mbox{for all $j=0,1,2$, \quad
and for all $x\in \mathcal O$.}\end{equation}

Simple computation yields
\[
{\phi}_{\alpha}=\frac{\partial \phi}{\partial
x_{\alpha}}={\phi}^{ij}v_{ij\alpha}, \ \
{\phi}_{\alpha\beta}=\frac{\partial^2 \phi}{\partial
x_{\alpha}\partial
x_{\beta}}={\phi}^{ij}v_{ij\alpha\beta}+{\phi}^{ij,km}v_{ij\alpha}v_{km\beta}.
\]
We differentiate equation (\ref{eq-v}) in $x_i$, by
(\ref{R-est1}),
\begin{equation}\label{orth1}
\sum_{\alpha \beta}F^{\alpha\beta}v_{\alpha\beta
i}+\sum_{k}F^{q_k}v_{k i}+F^{v}v_{ i}+ F^{x_i}=O(\phi),
\end{equation}
and differentiate equation (\ref{eq-v}) twice with respect to the
variables $x_i$ and $x_j$, again by (\ref{R-est1}),
\begin{eqnarray}\label{orth1-2}
&& \sum_{\alpha \beta} F^{\alpha\beta}v_{\alpha\beta
ij}+\sum_{\alpha\beta}v_{\alpha\beta i}
(\sum_{\gamma\eta}F^{\alpha\beta,\gamma\eta}v_{\gamma\eta j}+
\sum_kF^{\alpha\beta,q_k}v_{k j}+F^{\alpha\beta,v}v_{j}+
 F^{\alpha\beta,x_j} ) \nonumber
\\
&& +\sum_kF^{q_k}v_{k ij}+\sum_{k \alpha \beta}v_{k
i}(\sum_{\alpha \beta} F^{q_k,\alpha\beta}v_{\alpha\beta j}+
\sum_{l}F^{q_k,q_l}v_{lj}+F^{q_k,v}v_{j}+ F^{q_k,x_j}) \nonumber
\\
&& +F^{v}v_{ij}+v_{i}(\sum_{\alpha
\beta}F^{v,\alpha\beta}v_{\alpha\beta j}+
\sum_lF^{v,q_l}u_{lj}+F^{v,v}v_{j}+
 F^{v,x_j}) \nonumber
\\
&& +\sum_{\alpha \beta}F^{x_i,\alpha\beta}v_{\alpha\beta j}+
\sum_kF^{x_i,q_k}v_{kj}+F^{x_i,v}v_{j}+
 F^{x_i,x_j}=O(\phi).
\end{eqnarray}
As $v_{\alpha\beta ij}=v_{ij \alpha\beta}$ (this fact will have to
be modified later by a commutator formula when we deal with
symmetric curvature tensors on general manifolds), we get
\begin{eqnarray}\label{2-10}
\sum F^{\alpha\beta}\phi_{\alpha\beta}&=&\sum
F^{\alpha\beta}\phi^{ij}v_{ij\alpha\beta} +\sum
F^{\alpha\beta}\phi^{ij,km}v_{ij\alpha}v_{km\beta} \nonumber
\\
&=&\sum
F^{\alpha\beta}\phi^{ij,km}v_{ij\alpha}v_{km\beta}-\sum\phi^{ij}F^{q_k}v_{k
ij} \nonumber
\\
&& -\sum\phi^{ij}[F^v v_{ij}+2\sum
F^{\alpha\beta,q_k}v_{\alpha\beta i}v_{k j}+\sum F^{q_k,q_l}v_{k
i}v_{lj}\nonumber \\
&& +2\sum F^{q_k,v}v_{k i}v_{j}+ 2\sum F^{q_k,x_j}v_{k i}]
\nonumber
\\
&&-\sum \phi^{ij}[F^{\alpha\beta,\gamma\eta}v_{\alpha\beta
i}v_{\gamma\eta j}+2\sum F^{\alpha\beta,v}v_{\alpha\beta i}v_{j}+
2\sum F^{\alpha\beta,x_j}v_{\alpha\beta i} \nonumber \\
&& +\sum F^{v,v}v_{i}v_{j}+ 2 \sum F^{v,x_j}v_{j}+\sum F^{x_i
x_j}]+O(\phi)
\end{eqnarray}

We will deal terms in the right hand side of (\ref{2-10}). The
basic idea is to regroup them according indices in $G$ and $B$.
The analysis will be devoted to those third order derivatives
terms which have with at least two indices in $B$. Since it
contains some linear terms of such third order derivatives,
previous arguments in \cite{CGM} are not suitable here. The
introduction of function $q$ in (\ref{def-q}) is the key, the
concavity results of $q$ in last section will be used in crucial
way. As for the rest terms left in (\ref{2-10}), we will sort them
out in a way such that condition (\ref{cond-c}) can be used to
obtain appropriate control.

We note that since
$W=(v_{ij})$ is diagonal at $z$, by Lemma \ref{Lemma1-1} and Lemma
\ref{Lemma1-2},
\begin{equation}\label{2-11}
\phi^{ij}(z)=\left \{
\begin{array}{ll}
\sigma_l(G)+\frac{\sigma^2_{1}(B|i)-\sigma_{2}(B|i)}{\sigma^2_{1}(B)}+O(\phi),
&   {\rm if\ }i=j\in B \\
O(\phi),    &   {\rm if\ }i=j\in G   \\
0,& {\rm if\ } i\neq j
     \end{array}
   \right.
\end{equation}
Hence at $z$
\begin{eqnarray}\label{2-12}
&&\sum_{i,j}\phi^{ij}[F^v v_{ij}+2\sum
F^{\alpha\beta,q_k}v_{\alpha\beta i}v_{k j}+\sum F^{q_k,q_l}v_{k
i}v_{lj}+2\sum (F^{q_k,v}v_{k i}v_{j}+ F^{q_k,x_j}v_{k i})]\nonumber \\
&& =\sum_{i=1}^n\phi^{ii}[F^v v_{ii}+2\sum
F^{\alpha\beta,q_i}v_{\alpha\beta i}v_{i i}+ F^{q_i,q_i}v_{i
i}v_{ii}+2F^{q_i,v}v_{i i}v_{i}+
2F^{q_i,x_i}v_{ii}] \nonumber \\
&& =O(\phi)+\sum_{i\in B}\phi^{ii}[F^v +2\sum
F^{\alpha\beta,q_i}v_{\alpha\beta i}+ F^{q_i,q_i}v_{i
i}+2F^{q_i,v}v_{i}+ 2F^{q_i,x_i}]v_{ii} \nonumber \\
&   & \quad \quad  \leq O(\phi)+ C\sum_{i\in B}
(\sigma_l(G)+\frac{\sigma^2_{1}(B|i)-\sigma_{2}(B|i)}{\sigma^2_{1}(B)})v_{ii}=O(\phi).
\end{eqnarray}
This takes care of the third term in the right hand side of
(\ref{2-10}). For the second term there, we have
\begin{equation}\label{2-13}
\sum \phi^{ij}F^{q_k}v_{k ij}=O(\phi)+\sum_{i\in
B}\phi^{ii}F^{q_k}v_{k ii} =O(\phi+ \sum_{i,j\in B}|\nabla
v_{ij}|)
\end{equation}

For the fourth term in (\ref{2-10}), by (\ref{2-11}) we have,
\begin{eqnarray}\label{2-14}
&&\phi^{ij}[F^{\alpha\beta,\gamma\eta}v_{\alpha\beta
i}v_{\gamma\eta j}+2F^{\alpha\beta,v}v_{\alpha\beta i}v_{j}+
2F^{\alpha\beta,x_j}v_{\alpha\beta i}+F^{v,v}v_{i}v_{j}+ 2
F^{v,x_j}v_{j}+F^{x_i x_j}]\nonumber \\
&&=O(\phi)+ \sum_{i\in B}\phi^{ii} [\sum
F^{\alpha\beta,\gamma\eta}v_{\alpha\beta i}v_{\gamma\eta i}+2\sum
F^{\alpha\beta,v}v_{\alpha\beta i}v_{i}\nonumber \\
&& \quad + 2\sum F^{\alpha\beta,x_i}v_{\alpha\beta
i}+F^{v,v}v_{i}^2+ 2
F^{v,x_i}v_{i}+F^{x_i x_i}]\nonumber \\
&& =O(\phi+\sum_{i,j\in B}|\nabla v_{ij}|)+\sum_{i\in B}
(\sigma_l(G)+\frac{\sigma^2_{1}(B|i)-\sigma_{2}(B|i)}{\sigma^2_{1}(B)})\nonumber \\
&& \quad [\sum_{\alpha,\beta,\gamma,\eta\in G}
F^{\alpha\beta,\gamma\eta}v_{i\alpha\beta}v_{i\gamma\eta} +2
\sum_{\alpha,\beta\in G} F^{\alpha\beta,v}v_{i\alpha\beta}v_{i}+
2\sum_{\alpha,\beta\in G}F^{\alpha\beta,x_i}v_{i\alpha\beta}\nonumber \\
&& \quad +F^{v,v}v_{i}^2+ 2 F^{v,x_i}v_{i}+F^{x_i x_i}].
\end{eqnarray}

Now we deal with the term $\sum
F^{\alpha\beta}\phi^{ij,km}v_{ij\alpha}v_{km\beta}$ in
(\ref{2-10}). We note that
\[\phi^{ij,km}=\sigma_{l+1}^{ij,km}+q^{ij,km}.\] Since
$\sigma_{l-1}(W|ij)=O(\phi)$ for $i,j\in G,i\neq j$, for $\alpha,
\beta$ fixed, by (\ref{1-3}),
\begin{eqnarray*}
\sum \sigma_{l+1}^{ij,km}v_{ij\alpha}v_{km\beta}
 &=& \sum_{i\neq
k}\sigma_{l+1}^{ii,kk}v_{ii\alpha}v_{kk\beta} + \sum_{i\neq
j}\sigma_{l+1}^{ij,ji}v_{ij\alpha}v_{ji\beta}\\
& =& \sum_{i\neq k}\sigma_{l-1}(W|ik)v_{ii\alpha}v_{kk\beta}-
\sum_{i\neq j}\sigma_{l-1}(W|ij)v_{ij\alpha}v_{ji\beta}\\
&=& O(\phi+\sum_{i,j\in B}|\nabla v_{ij}|) - 2\sum_{i\in B,j\in
G}\sigma_{l-1}(G|j) v_{ij\alpha}v_{ij\beta}.
\end{eqnarray*}
As $\sigma_{l-1}(G|j)=\frac{\sigma_{l}(G)}{\lambda_j},j\in G$, we
have
\[
\sigma_{l+1}^{ij,km}v_{ij\alpha}v_{km\beta}=O(\phi+\sum_{i,j\in
B}|\nabla v_{ij}|) -2\sigma_{l}(G) \sum_{i\in B,j\in G}
\frac{1}{\lambda_j}v_{ij\alpha}v_{ij\beta}.
\]
By Proposition \ref{Prop1-2},
\begin{eqnarray*}
&& \sum_{i,j,k,m} q^{ij,km}v_{ij\alpha}v_{km\beta}
=O(\phi+\sum_{i,j\in B}|\nabla v_{ij}|) -2\sum_{i\in B,j\in G}
\frac{\sigma^2_{1}(B|i)-\sigma_{2}(B|i)}
{\sigma^2_{1}(B)\lambda_{j}}v_{ij\alpha}v_{ji\beta}\\
&& \quad  -\frac{1}{\sigma^3_{1}(B)}\sum_{i\in B}
(\sigma_{1}(B)v_{ii\alpha}-v_{ii}\sum_{j\in
B}v_{jj\alpha})(\sigma_{1}(B)v_{ii\beta}-v_{ii}\sum_{j\in
B}v_{jj\beta})
\\ && \quad - \frac{1}{\sigma_{1}(B)}\sum_{i,j\in B, i\neq j}
v_{ij\alpha}v_{ji\beta} -\frac{2}{\sigma^3_{1}(B)}\sum_{i \in B}
v_{ii}\sigma_{1}(B|i)v_{ii\alpha}v_{ii\beta}.
\end{eqnarray*}
We conclude that
\begin{eqnarray}\label{2-15}
&& \sum
F^{\alpha\beta}\phi^{ij,km}v_{ij\alpha}v_{km\beta}=O(\phi+\sum_{i,j\in
B}|\nabla
v_{ij}|)-\sum_{\alpha,\beta}F^{\alpha\beta}[\frac{2\sum_{i \in B}
v_{ii}\sigma_{1}(B|i)v_{ii\alpha}v_{ii\beta}}
{\sigma^3_{1}(B)}\nonumber
\\
&& \quad -\frac{1}{\sigma_{1}(B)} \sum_{i,j\in B, i\neq j}
v_{ij\alpha}v_{ji\beta} -2\sum_{i\in B}
(\sigma_l(G)+\frac{\sigma^2_{1}(B|i)-\sigma_{2}(B|i)}{\sigma^2_{1}(B)})
\frac{1}{\lambda_{j}}v_{ij\alpha}v_{ji\beta}\nonumber
\\
&& \quad  -\frac{1}{\sigma^3_{1}(B)}\sum_{i\in
B}(\sigma_{1}(B)v_{ii\alpha}-v_{ii}\sum_{j\in B}v_{jj\alpha})
(\sigma_{1}(B)v_{ii\beta}-v_{ii}\sum_{j\in
B}v_{jj\beta})].
\end{eqnarray}

Combining (\ref{2-12})-(\ref{2-15}), (\ref{2-10}) is deduced to
\begin{eqnarray}\label{new3-1}
&& \sum F^{\alpha\beta}\phi_{\alpha\beta}=O(\phi +\sum_{i,j\in
B}|\nabla v_{ij}|)-\frac{1}{\sigma_{1}(B)} \sum_{\alpha,\beta}
\sum_{i,j\in
B,i\neq j}F^{\alpha\beta}v_{ij\alpha}v_{ij\beta}\nonumber \\
&& \quad -\frac{2}{\sigma^3_{1}(B)}\sum_{\alpha,\beta}\sum_{i \in
B}
F^{\alpha\beta}v_{ii}\sigma_{1}(B|i)v_{ii\alpha}v_{ii\beta}\nonumber
\\
&& \quad -\frac{1}{\sigma_{1}^{3}(B)}\sum_{\alpha,\beta}\sum_{i\in
B} F^{\alpha\beta} (v_{ii\alpha}\sigma_{1}(B)-v_{ii}\sum_{j\in
B}v_{jj\alpha}) (v_{ii\beta}\sigma_{1}(B)-v_{ii}\sum_{j\in
B}v_{jj\beta})\nonumber \\
&& \quad  -\sum_{i\in B}[\sigma_l(G)+\frac{\sigma_{1}^{2}(B|i)-
\sigma_{2}(B|i)}{\sigma_{1}^{2}(B)}]
[\sum_{\alpha,\beta,\gamma,\eta\in
G}F^{\alpha\beta,\gamma\eta}(\Lambda)v_{i\alpha\beta}v_{i\gamma\eta}\nonumber
\\
&& \quad + 2\sum_{\alpha\beta\in G}F^{\alpha\beta} \sum_{j\in
G}\frac{1}{\lambda_{j}}v_{ij\alpha}v_{ij\beta}+2\sum_{\alpha,\beta\in
G}F^{\alpha\beta,v}v_{i\alpha\beta}v_{i} \nonumber \\
&& \quad + 2\sum_{\alpha,\beta\in
G}F^{\alpha\beta,x_i}v_{i\alpha\beta }+F^{v,v}v^2_{i}+ 2
F^{v,x_i}v_{i}+F^{x_i,x_i}].
\end{eqnarray}

At this point, we have succeeded in regrouping of terms involving
third order derivatives. We first estimate the fifth term on the
right hand side of (\ref{new3-1}). For each $i\in B$, let
\begin{eqnarray}\label{Ji}
J_i&=&[\sum_{\alpha,\beta,\gamma,\eta\in
G}F^{\alpha\beta,\gamma\eta}v_{i\alpha\beta}v_{i\gamma\eta}+
2\sum_{\alpha,\beta\in G}F^{\alpha\beta} \sum_{j\in
G}\frac{1}{\lambda_{j}}v_{ij\alpha}v_{ij\beta} \nonumber \\
&& \quad +2\sum_{\alpha,\beta\in G}F^{\alpha\beta,v}v_{
i\alpha\beta}v_{i}+ 2\sum_{\alpha,\beta\in
G}F^{\alpha\beta,x_i}v_{i\alpha\beta}+F^{v,v}v^2_{i}+ 2
F^{v,x_i}v_{i}+F^{x_i,x_i}].
\end{eqnarray}
If $l=0$, then $G=\emptyset$ and
\[
J_i=F^{v,v}(\nabla^2 v,\nabla v,v,z)v^2_{i}+ 2 F^{v,x_i}(\nabla^2
v,\nabla v,v,z)v_{i}+F^{x_i,x_i}(\nabla^2 v,\nabla v,v,z).
\]
Since $F\in C^{2,1}$ and $|\nabla^2 v(z)|=O(\phi)$, by condition
(\ref{wwcond-c}),
\[
J_i=F^{v,v}(0,\nabla v,v,z)v^2_{i}+ 2 F^{v,x_i}(0,\nabla
v,v,z)v_{i}+F^{x_i,x_i}(0,\nabla v,v,z)+O(\phi)\geq -C\phi.
\]
We may assume $1\le l\le n-1$. By Condition (\ref{cond-e}), since
$v\in C^{3,1}$ so $F^{\alpha \beta} \in C^{0,1}$, as
$\bar{\mathcal{O}} \subset \Omega$, there exists a constant
$\delta_{0}>0$, such that
\begin{equation}\label{cond-ee}
(F^{\alpha\beta})\geq \delta_{0}I,~~\forall y\in {\mathcal O}.
\end{equation}
As $l\geq 1$, so $n\in G$ and $F^{nn}\geq \delta_{0}$. From
(\ref{orth1}), since $v_{ik}=\delta_{ik}\lambda_i$ at $z$, we have
for $i\in B$
\[
\sum_{\alpha, \beta\in G}F^{\alpha\beta}v_{\alpha\beta i}+F^{v}v_{
i}+ F^{x_i}=O(\phi+\sum_{i,j\in B}|\nabla v_{ij}|),
\]
Now let's set $X_{\alpha\beta}=0$, $\alpha\in B$ or $\beta\in B$,
\[
X_{nn}=v_{inn}-\frac{1}{F^{nn}}[\sum_{\alpha, \beta\in
G}F^{\alpha\beta}v_{\alpha\beta i}+F^{v}v_{ i}+ F^{x_i}],
\]
$X_{\alpha\beta}=v_{i\alpha\beta}$ otherwise, $Y=-v_i$ and
$Z_k=-\delta_{ki}$. As $l\le n-1$, so that $(X_{\alpha\beta})\in
\mathcal{S}_{n-1}(\mbox{identity matrix})$ and
$\tilde{X}=((X_{\alpha\beta}), Y, Z_1,\cdots,Z_n) \in
\Gamma^{\bot}_{ X^*_F}$. Again by condition (\ref{wwcond-c}), we
infer that
\[J_i\geq -C(\phi+\sum_{i,j\in B}|\nabla v_{ij}|).\]

Since $C\ge \sigma_l(G)+\frac{\sigma_{1}^{2}(B|i)-
 \sigma_{2}(B|i)}{\sigma_{1}^{2}(B)}\geq 0$, thus we obtain
\begin{eqnarray}\label{2-16}
\sum_{\alpha,\beta} F^{\alpha\beta}\phi_{\alpha\beta}&\leq & C(
\phi +\sum_{i,j\in
B}|\nabla v_{ij}|) \nonumber \\
&-&\frac{1}{\sigma_{1}^{3}(B)}\sum_{\alpha,\beta}\sum_{i\in B}
F^{\alpha\beta}
 (v_{ii\alpha}\sigma_{1}(B)-v_{ii}\sum_{j\in B}v_{jj\alpha})
  (v_{ii\beta}\sigma_{1}(B)-v_{ii}\sum_{j\in B}v_{jj\beta})
\nonumber \\
&-&\frac{1}{\sigma_{1}(B)} \sum_{\alpha,\beta} \sum_{i,j\in
B,i\neq j}F^{\alpha\beta}v_{ij\alpha}v_{ij\beta}
-\frac{2}{\sigma^3_{1}(B)}\sum_{\alpha,\beta}\sum_{i \in B}
F^{\alpha\beta}v_{ii}\sigma_{1}(B|i)v_{ii\alpha}v_{ii\beta}.
\end{eqnarray}

The final stage of the proof is to control the term $\sum_{i,j\in
B}|\nabla v_{ij}|$ in (\ref{2-16}) by the rest terms on the right
hand side.
Let's set
\[
 V_{i\alpha}=v_{ii\alpha}\sigma_{1}(B)-v_{ii}\Big(\displaystyle\sum_{j\in
 B}v_{jj\alpha}\Big).
\]
By (\ref{cond-ee}),
\[
\sum_{\alpha,\beta} F^{\alpha\beta}V_{i\alpha}V_{i\beta}\geq
\delta_{0}\sum_{\alpha=1}^{n}V_{i\alpha}^{2}, \quad
\sum_{\alpha,\beta} F^{\alpha\beta}v_{ij\alpha}v_{ij\beta}\geq
\delta_{0}\sum_{\alpha=1}^{n}v_{ij\alpha}^{2}.
\]

Inserting above inequalities into (\ref{2-16}), we then obtain
\begin{eqnarray}\label{2-17}
\sum_{\alpha,\beta} F^{\alpha\beta}\phi_{\alpha\beta} &\leq & C(
\phi +\sum_{i,j\in B}|\nabla v_{ij}|)
-\frac{\delta_{0}}{\sigma_{1}^{3}(B)} \sum_{\alpha=1}^{n}\sum_{i\in
B} V_{i\alpha}^{2}\nonumber
\\
&& \quad -\frac{\delta_{0}}{\sigma_{1}(B)}
 \sum_{\alpha=1}^{n}\sum_{i,j\in B~ i\neq j} |v_{ij\alpha}|^{2}
 -\frac{2\delta_{0}}{\sigma^3_{1}(B)}\sum_{\alpha=1}^{n}\sum_{i \in B}
v_{ii}\sigma_{1}(B|i)v_{ii\alpha}^2.
\end{eqnarray}

The key differential inequality (\ref{2-6}) is the consequence of (\ref{2-17})
and the following lemma.

\begin{lemma}\label{le-control} There is a constant $C$ depending only on $n, \|v\|_{C^2}$
and $\frac{1}{\sigma_l(G)}$, such that
for any constant $ D>0$
\begin{equation}\label{control1}
\sum_{i,j\in B}|\nabla v_{ij}|\le C(1+\frac{2}{\delta_0}+D) (\phi+|\nabla \phi|)
+\sum_{\alpha=1}^{n}[\frac{\delta_{0}}{2}\frac{\sum_{i,j\in B~ i\neq j}|v_{ij\alpha}|^{2}}
{\sigma_{1}(B)}+ \frac{C}{D}\frac{\sum_{i\in
B}V_{i\alpha}^{2}}{\sigma_{1}^{3}(B)}].\end{equation}
\end{lemma}

\noindent{\bf Proof of Lemma \ref{le-control}.}\  We will use a
trick devised in \cite{Guan-Duke}. We break write
\[
\sum_{i,j\in B}|\nabla v_{ij}|=\sum_{i,j\in B,\ i\neq j}|\nabla
v_{ij}|+\sum_{i\in B}|\nabla v_{ii}|
\]
If $i\neq j$, for any $A>0$, the Cauchy-Schwarz inequality yields
\begin{eqnarray}
|v_{ij\alpha}|\leq
2\delta_{0}^{-1}\sigma_{1}(B)+\frac{\delta_{0}}{2}\frac{|v_{ij\alpha}|^{2}}
{\sigma_{1}(B)} \leq C\frac{2}{\delta_{0}}\phi +
\frac{\delta_{0}}{2}\frac{|v_{ij\alpha}|^{2}}
{\sigma_{1}(B)}.\nonumber
\end{eqnarray}

What left are the linear terms involving $v_{ii\alpha},~i\in B$,
we need the help of the second term on right hand side of
(\ref{2-17}) and $\phi_{\alpha}$. It follows from Lemma
\ref{Lemma1-2} that
\begin{equation}\label{2-18}
\phi_{\alpha}=O(\phi)+ \sum_{i\in B}(\sigma_{l}(G)
+\frac{\sigma_{1}^{2}(B|i)-
 \sigma_{2}(B|i)}{\sigma_{1}^{2}(B)})v_{ii\alpha}.
\end{equation}

Let us now fix $\alpha \in \{1,2,\cdots,n\}$, set
\[
P=\{i\in B|\ v_{ii\alpha}>0 \},\ N=\{i\in B|\ v_{ii\alpha}<0 \}, \
R=\{i\in B|\ v_{ii\alpha}=0 \}.
\] We consider two
separate cases.

{\bf Case 1.} Either $P=\emptyset$ or $N=\emptyset$. In this case,
$v_{ii\alpha}$ has the same sign for all $i\in B$. We can derive
easily
\begin{equation}
|v_{ii\alpha}|= O(\phi + |\phi_{\alpha}|).
\end{equation}

{\bf Case 2.} $P\neq \emptyset,~N\neq\emptyset$. We may assume
\[
\sum_{i\in P} v_{ii}\geq\sum_{j\in N} v_{jj},
\]
by reversing the direction of $\partial_{x_{\alpha}}$ if
necessary, since we only need to control $|v_{ii\alpha}|$. It
follows from (\ref{2-18}) that, for $i\in P$,
 $$v_{ii\alpha}\leq \sum_{k\in
 P}v_{kk\alpha} \leq
 \frac{1}{\sigma_{l}(G)}O(\phi+|\phi_{\alpha}|)
 -C\sum_{j\in N}v_{jj\alpha},$$
 for some positive constant $C$ under control. At this point, we have switched the
 estimation of $v_{ii\alpha},~i\in P$ to the estimation of $-v_{jj\alpha},~j\in
 N$.

\vskip0.4cm \noindent {\bf Claim:} ~~If
$P\neq\emptyset,~N\neq\emptyset,~\sum_{i\in P} v_{ii}\geq\sum_{j\in
N} v_{jj}$, we have
\[
 \Big(\sum_{j\in
N}v_{jj\alpha}\Big)^{2} \leq
\frac{4n^{2}}{\sigma_{1}^{2}(B)}\sum_{i\in B}V_{i\alpha}^{2}.
\]

If the {\bf Claim} is true, we get for all $k\in N$,
\begin{eqnarray}\label{new2-2}
-v_{kk\alpha}&\leq &-\sum_{j\in N}v_{jj\alpha} \nonumber \\
&\leq & D \sigma_{1}(B)+ \frac{\Big(\sum_{j\in
N}v_{jj\alpha}\Big)^{2}}{D
\sigma_{1}(B)}\nonumber \\
&\leq & CD\phi +
\frac{4n^2}{D}\frac{1}{\sigma_{1}^{3}(B)}\sum_{i\in
B}V_{i\alpha}^{2}.
\end{eqnarray}
which can be controlled by the 3rd term in (\ref{2-17}) if we choose
the constant $D$ large enough. Consequently we can control terms
involving $v_{ii\alpha},~i\in P$.  We now validate the {\bf Claim}.

\vskip0.2cm
 \noindent {\it Proof of {\bf Claim}.} \ We first have by the Cauchy-Schwarz inequality
\[
\Big(\sum_{i\in N}V_{i\alpha}\Big)^{2}\leq n^2 \sum_{i\in
N}V_{i\alpha}^{2}\leq n^2\sum_{i\in B}V_{i\alpha}^{2}.
\]
It follows that from the definitions of the sets $P,N,R$ and
$V_{i\alpha}$
\begin{eqnarray}\label{neweq1}
-\sum_{i\in N}V_{i\alpha}&=&\sum_{i\in N}\Big(v_{ii}(\sum_{j\in
N}v_{jj\alpha}+\sum_{k\in P}v_{kk\alpha})-v_{ii\alpha}(\sum_{j\in
N}v_{jj}+\sum_{j\in R}v_{jj}+\sum_{k\in P}v_{kk})\Big) \nonumber \\
&&=\Big(\sum_{i\in N}v_{ii}\Big)\Big(\sum_{k\in
P}v_{kk\alpha}\Big)-\Big(\sum_{k\in P\cup
R}v_{kk}\Big)\Big(\sum_{i\in N} v_{ii\alpha}\Big)
\end{eqnarray}
Since in this case
\[
\sum_{i\in N}v_{ii}\geq 0, \sum_{k\in P}v_{kk\alpha}>0, \sum_{j\in
N}v_{jj\alpha}\leq 0,
\]
all the terms on the right hand side of (\ref{neweq1}) are
nonnegative, thus we obtain
\[
\Big(\sum_{i\in N}V_{i\alpha}\Big)^{2}\geq \Big(\sum_{k\in P\cup
R}v_{kk}\Big)^{2}\Big(\sum_{i\in N}v_{ii\alpha}\Big)^{2} \geq
\Big(\frac{1}{2}\sum_{k\in B}v_{kk}\Big)^{2}\Big(\sum_{i\in
N}v_{ii\alpha}\Big)^{2}= \frac{\sigma_{1}^{2}(B)}{4}\Big(\sum_{i\in
N}v_{ii\alpha}\Big)^{2}.
\]
The lemma is proved.\qed
\medskip

By Lemma \ref{le-control} and (\ref{2-17}), there exist positive constants
$C_1, C_2$ independent of
$\epsilon$, such that
\begin{eqnarray}\label{2-17nn}
\sum_{\alpha,\beta} F^{\alpha\beta}\phi_{\alpha\beta} \leq  C_1(
\phi +|\nabla \phi|)-C_2\sum_{i,j \in B} |\nabla v_{ij}|.
\end{eqnarray}
Taking $\epsilon \to 0$, (\ref{2-17nn}) is proved for $u$. By the
Strong Maximum Principle, $\phi \equiv 0$ in $\mathcal O$. Since
$\Omega$ is flat, following the arguments in \cite{CF85, Kl87},
for any $x_0\in \Omega$, there is a neighborhood $\mathcal U$ and
$(n-l)$ fixed directions $V_1, \cdots, V_{n-l}$ such that
$\nabla^2 u(x)V_j=0$ for all $1\le j\le n-l$ and $x\in \mathcal
U$. The proof of Theorem \ref{microc-w} is complete. \qed

\section{Condition (\ref{cond-c}) and discussions}
We discuss the convexity condition (\ref{cond-c}) in this section.
We write $A^{-1}=(A^{ij})$ to be the inverse matrix $A^{-1}$ of
positive definite matrix $A$.

\begin{lemma}\label{Lemma2-1}\ $F$ satisfies Condition (\ref{cond-c}) if and only if
\begin{eqnarray}\label{eq2-2}
\quad & & \sum_{i,j,k,l=1}^n F^{ij,kl}(A,p,u,x)X_{ij}X_{kl}
+2\sum_{i,j,k,l=1}^n F^{ij}(A,p,u,x)A^{kl} X_{ik}X_{jl}+F^{u,u}Y^2
\nonumber \\
& & -2\sum_{i,j=1}^n F^{ij,u}X_{ij}Y-2\sum_{i,j,k=1}^n
F^{ij,x_k}X_{ij}Z_k + 2\sum_{i=1}^n F^{u,x_i}YZ_i+ \sum_{i,j=1}^n
F^{x_i,x_j}Z_iZ_j\geq 0
\end{eqnarray}
for every $X=(X_{ij})\in \mathcal{S}^n$, $Y\in \mathbb R$ and
$Z=(Z_i)\in\mathbb R^n$.
\end{lemma}

\noindent {\bf Proof.}\ We have, from the convexity of $\tilde
F(B,u,x)=F(B^{-1},u, p, x)$ (for each $p$ fixed),
\begin{eqnarray}\label{eq2-3}
\quad & & \sum_{\alpha,\beta,\gamma,\eta=1}^n
\tilde{F}^{\alpha\beta,\gamma\eta}(B,u,x)\tilde{X}_{\alpha\beta}\tilde{X}_{\gamma\eta}
+2\sum_{\alpha,\beta=1}^n\tilde{F}^{\alpha\beta,u}\tilde{X}_{\alpha\beta}Y+\tilde{F}^{u,u}Y^2
\nonumber \\
& &
+2\sum_{\alpha,\beta,k=1}^n\tilde{F}^{\alpha\beta,x_k}\tilde{X}_{\alpha\beta}Z_k
+ 2\sum_{k=1}^n
\tilde{F}^{u,x_k}YZ_k+\sum_{i,j=1}^nF^{x_i,x_j}Z_iZ_j \geq 0
\end{eqnarray}
for every $\tilde{X}\in \mathcal{S}^{n}$, $Y\in \mathbb R$,
$Z=(Z_i)\in\mathbb R^n$ and $B\in \mathcal{S}^{n}_{+}$. A direct
computation yields
\[
\tilde{F}^{\alpha\beta}(B,u,x)=-F^{ij}(B^{-1},p,u,x)B^{i\alpha}
B^{j\beta},
\]
\[
\tilde{F}^{\alpha\beta,u}(B,u,x)=-F^{ij,u}(B^{-1},p,u,x)B^{i\alpha}
B^{j\beta},
\]
\[
\tilde{F}^{\alpha\beta,\gamma\eta}(B,u,x)=F^{ij,kl}(B^{-1},p,u,x)B^{i\alpha}
B^{j\beta}B^{k\gamma} B^{l\eta}
\]
\[
+F^{ij}(B^{-1},p,u,x)(B^{i\gamma}B^{j\beta}B^{\eta\alpha}
+B^{i\alpha}B^{j\eta}B^{\beta\gamma}).
\]
Other derivatives can be calculated in a similar way. Substituting
these into (\ref{eq2-3}), (\ref{eq2-2}) follows directly. \qed

\medskip

Let $Q \in \mathbb O_{n}$, we define
\[ \tilde{F}_Q(A,u,x)=F(Q\left (
\begin{array}{cc}
       0 & 0 \\ 0 & A^{-1}  \end{array} \right )Q^T,p,u,x)
\]
for $(A,u,x)\in \mathcal{S}^{n-1}_+ \times \mathbb{R}\times \Omega$
and fixed $p$. Condition (\ref{cond-c}) implies the following
condition
\begin{eqnarray}\label{wcond-c}
\tilde{F}_Q(A,u,x) \quad \mbox{is locally convex}
\end{eqnarray}
in $\mathcal{S}^{n-1}_+\times \mathbb{R}\times \Omega$ for any
fixed $n\times n$ orthogonal matrix $Q$.

Lemma \ref{Lemma2-1} yields the following by approximating.

\begin{corollary}\label{corollary2-2}\ Let $Q \in \mathbb O_{n}$. Assume $F$ satisfies
condition (\ref{wcond-c}), then
\begin{eqnarray}\label{eq2-4}
Q^*(\tilde X, \tilde X) \ge 0,
\end{eqnarray}
for every $\tilde X=((X_{ij}), Y, Z_1,\cdots,Z_n) \in
\mathcal{S}_{n-1}(Q)\times \mathbb R\times \mathbb R^n$, where
$Q^{*}$ is defined in (\ref{eq2-5}).
\end{corollary}

\medskip

In particular, by Corollary \ref{corollary2-2}, condition
(\ref{wcond-c}) implies (\ref{wwcond-c}). Since condition
(\ref{cond-c}) implies (\ref{wcond-c}), Lemma \ref{lemma-cdc} is a
consequence of Corollary \ref{corollary2-2}.

Condition (\ref{wcond-c}) is weaker than condition (\ref{cond-c}).
In particular condition (\ref{wcond-c}) is empty condition in $A$
when $n=1$. There is a wide class of functions which satisfy
(\ref{eq2-4}). The most important examples are $\sigma_k$ and
$\frac{\sigma_l}{\sigma_k}$ ($l>k$). If $g$ is convex and $F_1,\cdots, F_m$ are in this
class, then $F=g(F_1,\cdots,F_m)$ is also in this class. In particular,
if $F_1>0$ and $F_2>0$ are in the class, so is
$F=F_1^{\alpha}+F_2^{\beta}$ for any $\alpha\ge 1$, $\beta\ge 1$.
Another property of condition (\ref{wcond-c}) is the following

\begin{corollary}
\label{cor1c} If $F$ satisfies (\ref{eq2-4}), then so is the
function $G(A)=F(A+E)$ for any nonnegative definite matrix $E$.
\end{corollary}

We also have the following lemma.

\begin{lemma}\label{homog2}
If $n=2$ and $F(A)\ge 0$ is symmetric and of homogeneous of degree
$k$. If either $k\le0$ or $k\ge 1$, then $F$ satisfies
(\ref{eq2-4}).
\end{lemma} \noindent {\bf Proof.} Since $n=2$, condition (\ref{eq2-4}) is equivalent to
$F^{\lambda_2, \lambda_2}\ge 0$. By homogeneity, we have
\[\sum_{i,j=1}^n F^{\lambda_i, \lambda_j}\lambda_i\lambda_j=k(k-1)F.\]
$n=2$ and $\lambda_1=0$ yields $F^{\lambda_2,
\lambda_2}\lambda_2^2=k(k-1)F(0,\lambda_2)\ge 0$. \qed

\bigskip

Simple example like $u=\sum_{i=1}^n x_i^4$, $F(A)=\sigma_1(A)$
indicates that certain condition is needed in Theorem
\ref{microc}. If $F$ is independent of $x,u$, one may ask if the
convexity assumption of $F(A^{-1},p)$ for $A$ in condition
(\ref{cond-c}) (or condition \ref{wwcond-c}) is necessary for
Theorem \ref{microc}. As we remarked before, when $n=1$, it is not
necessary. For general $n\ge 2$, we have the following theorem.

\begin{theorem}\label{n2}
Suppose $F(A,p)$ is elliptic and $u$ is a convex solution of
\begin{equation}\label{homo1}
F(\nabla^2 u, \nabla u)=0,\end{equation}
then $W=(\nabla^2 u)$ is either of
constant rank, or its minimal rank is at least $2$. In particular,
if $n=2$, then $W$ is of constant rank. \end{theorem}

\noindent{Proof.} The proof follows same lines of proof of Theorem
\ref{microc-w} with the following observations: condition (\ref{wcond-c})
was only used to control $J_i$ defined in (\ref{Ji}). Let $l$ be the
minimum rank of $W$. If $l=0$, that is $G=\emptyset$, the proof of
Theorem \ref{microc-w} works without any change since $F$ is independent of $(u,x)$
in our case. What left is the
case $l=1$, i.e., $|G|=1$, we may assume $\alpha=n\in G$. Note
that (\ref{new3-1}) still holds. Since $F(\nabla^2 u, \nabla
u)=0$, and
\[0=\nabla_i F(\nabla^2 u, \nabla u)=  F^{nn}u_{nni}+O(\phi+\sum_{i,j\in B}|\nabla u_{ij}|).\]
This gives
\[|u_{nni}|\le C(\phi +\sum_{i,j\in B}|\nabla u_{ij}|).\]
Of course, the treatment of terms involving $u_{ij\beta}$ for
$i,j\in B$ follows the same way as in the proof of Theorem
\ref{microc-w}. We can now deduce that $W$ is of constant.
Finally, if $n=2$, the only other case is $l=2$. In this case, $W$
is of full rank everywhere. \qed

\medskip
\begin{remark} In \cite{BL}, Bramscap and Lieb proved $\log$-concavity of the
first eigenfunction of Laplacian operator for bounded convex
domains in $\mathbb R^n$ (see also \cite{SWYY, CS} for different
proofs). In general, for a nonlinear eigenvalue problem
$F(\nabla^2 v)=\lambda v$, the function $u=-\log v$ satisfies
equation (\ref{homo1}) if $F$ is of homogeneous degree of one.
\end{remark}

\begin{remark}\label{rkn2} The above proof of Theorem \ref{n2}
indicates that if the minimal rank of $W$ is either $0$ or $1$,
then the rank of $(\nabla^2 u)$ is the same everywhere. There is
no structure condition imposed on $F$ except the ellipticity
condition (\ref{cond-e}). This observation will be used in the
proof of Theorem \ref{codazzi-1} in the next section. \end{remark}

\medskip

We conclude this section with the proof of Theorem
\ref{microc-para}. It is a consequence of the following Strong
Maximum Principle for parabolic equations.

\begin{theorem}\label{microc-w-para} Suppose that the function
$F\in C^{2,1}$ satisfies conditions (\ref{cond-e}) and
(\ref{eq2-4}) for each $t\in [0,T]$, let $u \in C^{3}(\Omega\times
[0, T])$ is a convex solution of (\ref{equ1-parab}). For each
$0<t_0\le T$, if $\nabla^2 u$ attains minimum rank $l$ at certain
point $x_0\in \Omega$, then there exist a neighborhood $\mathcal
O$ of $x_0$ and a positive constant $C$ independent of
$\phi$ (defined in (\ref{def-phi})), such that for $t$ close to
$t_0$, $\sigma_{l}(u_{ij}(x,t))>0$ for $x\in \mathcal O$, and
\begin{equation}\label{2-6-p}
\sum_{\alpha, \beta}
F^{\alpha\beta}\phi_{\alpha\beta}(x,t)-\phi_t(x,t) \leq
 C(\phi(x,t)+|\nabla\phi(x,t)|), \quad \forall x\in \mathcal O.
\end{equation}
Consequently, the rank of $\nabla^2 u(x,t)$ is constant for every
fixed $t>0$ and it is non-decreasing. For each $0<t\le T$, $x_0\in
\Omega$, there exist a neighborhood $\mathcal U$ of $x_0$ and
$(n-l(t))$ fixed directions $V_1, \cdots, V_{n-l(t)}$  such that
$\nabla^2u(x,t)V_j=0$ for all $1\le j\le n-l(t)$ and $x\in
\mathcal U$. Furthermore, for any $t_0$, there is $\delta>0$, such
that the null space of $\nabla^2 u(x,t)$ is parallel for $(x,t)\in
\mathcal O \times (t_0,t_0+\delta)$.
\end{theorem}

\noindent {\bf Proof of Theorem \ref{microc-w-para}.} The proof is
similar to the proof of Theorem \ref{microc-w}, here we will use
the Strong Maximum Principle for parabolic equations.

Since $u\in C^{3}$, and the assumption on $F$, $u\in C^4$
automatically. Suppose $(\nabla^2u(x,t_0))$ attains minimal rank
$l$ at some point $x_0\in \Omega$. We may assume $l\le n-1$,
otherwise there is nothing to prove. By continuity,
$\sigma_{l}(u_{ij}(x,t))>0$ in a neighborhood of $(x_0,t_0)$. We
want to show (\ref{2-6-p}).

With
$u_t=F(\nabla^2 u, \nabla u, u, x, t)$, using the same notations
as in the proof of Theorem \ref{microc-w}, equation
(\ref{orth1-2}) becomes
\begin{eqnarray}\label{orth1-2-t}
&& \sum_{\alpha \beta} F^{\alpha\beta}v_{\alpha\beta
ij}+\sum_{\alpha\beta}v_{\alpha\beta i}
(\sum_{\gamma\eta}F^{\alpha\beta,\gamma\eta}v_{\gamma\eta j}+
\sum_kF^{\alpha\beta,q_k}v_{k j}+F^{\alpha\beta,v}v_{j}+
 F^{\alpha\beta,x_j} ) \nonumber
\\
&& +\sum_kF^{q_k}v_{k ij}+\sum_{k \alpha \beta}v_{k
i}(\sum_{\alpha \beta} F^{q_k,\alpha\beta}v_{\alpha\beta j}+
\sum_{l}F^{q_k,q_l}v_{lj}+F^{q_k,v}v_{j}+ F^{q_k,x_j}) \nonumber
\\
&& +F^{v}v_{ij}+v_{i}(\sum_{\alpha
\beta}F^{v,\alpha\beta}v_{\alpha\beta j}+
\sum_lF^{v,q_l}u_{lj}+F^{v,v}v_{j}+
 F^{v,x_j}) \nonumber
\\
&& +\sum_{\alpha \beta}F^{x_i,\alpha\beta}v_{\alpha\beta j}+
\sum_kF^{x_i,q_k}v_{kj}+F^{x_i,v}v_{j}+
 F^{x_i,x_j}=O(\phi)+v_{ij,t},
\end{eqnarray}
and accordingly, equation (\ref{2-10}) becomes
\begin{eqnarray}\label{2-10-t}
\sum F^{\alpha\beta}\phi_{\alpha\beta}&=&\sum
F^{\alpha\beta}\phi^{ij}v_{ij\alpha\beta} +\sum
F^{\alpha\beta}\phi^{ij,km}v_{ij\alpha}v_{km\beta} \nonumber
\\
&=&\sum
F^{\alpha\beta}\phi^{ij,km}v_{ij\alpha}v_{km\beta}-\sum\phi^{ij}F^{q_k}v_{k
ij} \nonumber
\\
&& -\sum\phi^{ij}[F^v v_{ij}+2\sum
F^{\alpha\beta,q_k}v_{\alpha\beta i}v_{k j}+\sum F^{q_k,q_l}v_{k
i}v_{lj}\nonumber \\
&& +2\sum F^{q_k,v}v_{k i}v_{j}+ 2\sum F^{q_k,x_j}v_{k i}]
\nonumber
\\
&&-\sum \phi^{ij}[F^{\alpha\beta,\gamma\eta}v_{\alpha\beta
i}v_{\gamma\eta j}+2\sum F^{\alpha\beta,v}v_{\alpha\beta i}v_{j}+
2\sum F^{\alpha\beta,x_j}v_{\alpha\beta i} \nonumber \\
&& +\sum F^{v,v}v_{i}v_{j}+ 2 \sum F^{v,x_j}v_{j}+\sum F^{x_i
x_j}]+O(\phi)+\sum \phi^{ij}v_{ij,t}
\end{eqnarray}
We note that $\phi_t=\sum \phi^{ij}v_{ij,t}$, equation
(\ref{2-10-t}) can be written as
\begin{eqnarray}\label{2-10-tt}
\sum F^{\alpha\beta}\phi_{\alpha\beta}-\phi_t &=&\sum
F^{\alpha\beta}\phi^{ij,km}v_{ij\alpha}v_{km\beta}-\sum\phi^{ij}F^{q_k}v_{k
ij} \nonumber
\\
&& -\sum\phi^{ij}[F^v v_{ij}+2\sum
F^{\alpha\beta,q_k}v_{\alpha\beta i}v_{k j}+\sum F^{q_k,q_l}v_{k
i}v_{lj}\nonumber \\
&& +2\sum F^{q_k,v}v_{k i}v_{j}+ 2\sum F^{q_k,x_j}v_{k i}]
\nonumber
\\
&&-\sum \phi^{ij}[F^{\alpha\beta,\gamma\eta}v_{\alpha\beta
i}v_{\gamma\eta j}+2\sum F^{\alpha\beta,v}v_{\alpha\beta i}v_{j}+
2\sum F^{\alpha\beta,x_j}v_{\alpha\beta i} \nonumber \\
&& +\sum F^{v,v}v_{i}v_{j}+ 2 \sum F^{v,x_j}v_{j}+\sum F^{x_i
x_j}]+O(\phi)
\end{eqnarray}
Now the right hand side of (\ref{2-10-tt}) is the same as the
right hand side of (\ref{2-10}). The same analysis in the proof of
Theorem \ref{microc-w} for the right hand side of equation
(\ref{2-10}) yields
\begin{equation}\label{2-9-t}
 \sum F^{\alpha\beta}\phi_{\alpha\beta}(x,t)-\phi_t(x,t) \leq
 C_{1}(\phi(x,t)+|\nabla\phi(x,t)|)-C_2 \sum_{i,j \in B} |\nabla v_{ij}|.
\end{equation}
We now $\nabla^2u(x,t)$ is of constant rank $l(t)$ for each $t>0$.
Since w$\Omega$ is flat, by the arguments in \cite{CF85, Kl87},
for each $0<t\le T$, $x_0\in \Omega$, there exist a neighborhood
$\mathcal U$ of $x_0$ and $(n-l(t))$ fixed directions $V_1,
\cdots, V_{n-l(t)}$  such that $\nabla^2u(x,t)V_j=0$ for all $1\le
j\le n-l(t)$ and $x\in \mathcal U$. Now back to (\ref{2-9-t}), we
have $\sum_{i,j \in B} |\nabla u_{ij}(x,t)|\equiv 0$, therefore
the null space of $\nabla^2u$ is parallel.  \qed

\begin{remark}\label{wconds}
Tracing back to our proofs, for Theorem \ref{microc}, we only need
locally convexity condition in (\ref{cond-c}) near solution $u$ at
the points where some of eigenvalues of $\nabla^2 u$ are small.
For solution $u$ of (\ref{equ1}), we let
\begin{equation}\label{gamma}
\mathcal{D}_{u(x)} = \{
\mbox{$r$ diagonal} | \quad r=Q(\nabla^2u(x))Q^T  \quad \mbox{for some $Q\in
O(n)$}\}.\end{equation}
For each $\delta>0$, set $I_{u(x)}^{\delta}=\{s| \quad |s-u(x)|\le \delta\}$, and
\[\tilde D_{u(x)}^{\delta}=\{A| \quad \|A^{-1}-r\|\le \delta,
\mbox{ for some $r\in \mathcal{D}_{u(x)}$}\}.\] The condition
(\ref{cond-c}) in Theorem \ref{microc} can be replaced by: there
is $\delta>0$ and for $p=Q\nabla u(x)$ ($Q\in O(n)$),
\begin{equation}\label{cond-c-1} F(A^{-1},p,u,x) \quad \mbox{is
locally convex in $(A, u, x)$ in $\tilde{D}_{u(x)}^{\delta}\times
I_{u(x)}^{\delta}\times \mathcal O$}.
\end{equation}
Similarly, for condition (\ref{cond-cp}) and condition
(\ref{wcond-c}) are only needed to be valid for $(A, u, x)$ in
$\tilde{D}_{u(x)}^{\delta}\times I_{u(x)}^{\delta}\times
\mathcal O$ for each $t$. We also remark that regularity assumptions on $u$ and $F$ in
Theorem \ref{microc-para} and Theorem \ref{microc-w-para} can be
reduced to be $C^2$.
\end{remark}

\section{Geometric applications}

We discuss geometric nonlinear differential equations in this section.

\begin{proposition}\label{thmw2-pro}
Suppose $F(A,X, \vec n, t)$ is elliptic in $A$ and satisfies
condition (\ref{eq2-4}) for each fixed $\vec n\in \mathbb S^n$,
$t\in [0,T]$ for some $T>0$. Let $M(t)$ be oriented immersed
connect hypersurface in $ \mathbb R^{n+1}$ with a nonnegative
definite second fundamental form $h(t)$ satisfying equation
(\ref{flow1}), then $h(t)$ is of constant rank for each
$t\in(0,T]$. Moreover, if let $l(t)$ be the minimal rank of
$h(t)$, then $l(s)\le l(t)$ for all $0<s\le t\le T$ and the null
space of $h$ is parallel for each $t$.
\end{proposition}

We note that Theorem \ref{thmw2} follows directly from Proposition \ref{thmw2-pro}
(since equation (\ref{w3.10}) is a special case of equation (\ref{flow1}) by making
$M$ independent
of $t$) and a splitting theorem for complete hypersurface in $\mathbb R^{n+1}$.

\medskip

\noindent{\bf Proof of Proposition \ref{thmw2-pro}.}  For $\epsilon>0$, let
$W=(g^{im}h_{mj}+\epsilon \delta_{ij}) $, where $h=(h_{ij})$ the second
fundamental form and $(g_{ij})$ the first fundamental form of $M(t)$, and let $l(t)$ be
the minimal rank of $h(t)$. For a fixed $t_0\in (0,T)$, let $x_0\in M$ such that $h(t_0)$
attains minimal rank
at $x_0$. Set
$\phi(x,t)=\sigma_{l+1}(W(x,t))+\frac{\sigma_{l+2}}{\sigma_{l+1}}(W(x,t))$.
$\phi$ is in $C^{1,1}$ by result of section 2.
We want establish that in a small neighborhood of $(x_0,t_0)$, there are constants
$C_1, C_2$ independent of $\epsilon$ such that
\begin{equation}\label{flow3}
F^{ij}\phi_{ij}-\phi_t\le C_1\phi+C_2|\nabla \phi|.\end{equation}
The proposition follows from (\ref{flow3}) and the Strong Maximum Principle for
parabolic equations by taking $\epsilon\to 0$.

We work on $W=(h_{ij}+\epsilon g_{ij}) $ in place of Hessian
$(v_{ij})$ in the proof of Theorem \ref{microc-w}. We set position
vector $X=(X^1,\cdots, X^{n+1})$. (\ref{flow3}) can be proved
using the arguments in the proofs of Theorem \ref{microc-w} and
Theorem \ref{microc-para} and the Gauss equation, Codazzi equation
and the Weingarten equation for hypersurfaces. We note that under
(\ref{flow1}), the Weingarten form $h^i_j=g^{im}h_{mj}$ satisfies
equation
\begin{equation}\label{flow2} \partial_t h^i_j  =
\nabla^i\nabla_j F + F (h^2)^i_j,\end{equation} where
$h^2=(h^i_lh^l_j)$.

The same arguments in the proof of Theorem \ref{microc-w} can carry through
some modifications to get parabolic version of (\ref{orth1-2}) using (\ref{flow2}).
In this case, $W_{ijkm}$ and $W_{kmij}$ may be different.
But as $W$ is Codazzi, the commutator term can be controlled using
the Ricci identity. Also, $p$ is replaced by $\vec
n$, we use the Gauss equation when we differentiate in $p$. All
these terms are controlled by $CW_{ii}$. We notice that $W_{ii}\le
\phi$ for all $i\in B$, so we have the following corresponding
formula to replace (\ref{new3-1}),
\begin{eqnarray}\label{new3-1-n}
&& \sum F^{\alpha\beta}\phi_{\alpha\beta}-\phi_t=O(\phi +\sum_{i,j\in
B}|\nabla W_{ij}|)-\frac{1}{\sigma_{1}(B)} \sum_{\alpha,\beta}
\sum_{i,j\in
B,i\neq j}F^{\alpha\beta}W_{ij\alpha}W_{ij\beta}\nonumber \\
&& \quad -\frac{2}{\sigma^3_{1}(B)}\sum_{\alpha,\beta}\sum_{i \in
B}
F^{\alpha\beta}W_{ii}\sigma_{1}(B|i)W_{ii\alpha}W_{ii\beta}\nonumber
\\
&& \quad -\frac{1}{\sigma_{1}^{3}(B)}\sum_{\alpha,\beta}\sum_{i\in
B} F^{\alpha\beta} (W_{ii\alpha}\sigma_{1}(B)-W_{ii}\sum_{j\in
B}v_{jj\alpha}) (W_{ii\beta}\sigma_{1}(B)-W_{ii}\sum_{j\in
B}v_{jj\beta})\nonumber \\
&& \quad  -\sum_{i\in B}[\sigma_l(G)+\frac{\sigma_{1}^{2}(B|i)-
\sigma_{2}(B|i)}{\sigma_{1}^{2}(B)}]
[\sum_{\alpha,\beta,\gamma,\eta\in
G}F^{\alpha\beta,\gamma\eta}(\Lambda)W_{i\alpha\beta}W_{i\gamma\eta}
+\sum_{\alpha}F^{X^{\alpha}}X^{\alpha}_{ii}\nonumber
\\
&& \quad + 2\sum_{\alpha\beta\in G}F^{\alpha\beta} \sum_{j\in
G}\frac{1}{\lambda_{j}}W_{ij\alpha}W_{ij\beta}+
2\sum_{\alpha,\beta\in
G}\sum_{\gamma=1}^{n+1}F^{\alpha\beta,X^{\gamma}}W_{i\alpha\beta
}X^{\gamma}_i +\sum_{\gamma, \eta=1}^{n+1}
F^{X^{\gamma},X^{\eta}}X^{\gamma}_iX^{\eta}_i].
\end{eqnarray}

The term involving $X_{ii}$ is controlled by $Ch_{ii}$ (in turn by
$CW_{ii}$) using the Weingarten formula. We obtain

\begin{eqnarray}\label{new3-1-nn-f}
&& \sum F^{\alpha\beta}\phi_{\alpha\beta}-\phi_t=O(\phi +\sum_{i,j\in
B}|\nabla W_{ij}|)-\frac{1}{\sigma_{1}(B)} \sum_{\alpha,\beta}
\sum_{i,j\in
B,i\neq j}F^{\alpha\beta}W_{ij\alpha}W_{ij\beta}\nonumber \\
&& \quad -\frac{2}{\sigma^3_{1}(B)}\sum_{\alpha,\beta}\sum_{i \in
B}
F^{\alpha\beta}W_{ii}\sigma_{1}(B|i)W_{ii\alpha}W_{ii\beta}\nonumber
\\
&& \quad -\frac{1}{\sigma_{1}^{3}(B)}\sum_{\alpha,\beta}\sum_{i\in
B} F^{\alpha\beta} (W_{ii\alpha}\sigma_{1}(B)-W_{ii}\sum_{j\in
B}v_{jj\alpha}) (W_{ii\beta}\sigma_{1}(B)-W_{ii}\sum_{j\in
B}v_{jj\beta})\nonumber \\
&& \quad  -\sum_{i\in B}[\sigma_l(G)+\frac{\sigma_{1}^{2}(B|i)-
\sigma_{2}(B|i)}{\sigma_{1}^{2}(B)}]
[\sum_{\alpha,\beta,\gamma,\eta\in
G}F^{\alpha\beta,\gamma\eta}(\Lambda)W_{i\alpha\beta}W_{i\gamma\eta}
\nonumber
\\
&& \quad + 2\sum_{\alpha\beta\in G}F^{\alpha\beta} \sum_{j\in
G}\frac{1}{\lambda_{j}}W_{ij\alpha}W_{ij\beta}+
2\sum_{\alpha,\beta\in
G}\sum_{\gamma=1}^{n+1}F^{\alpha\beta,X^{\gamma}}W_{i\alpha\beta
}X^{\gamma}_i +\sum_{\gamma, \eta=1}^{n+1}
F^{X^{\gamma},X^{\eta}}X^{\gamma}_iX^{\eta}_i].
\end{eqnarray}
The right hand side of (\ref{new3-1-nn-f}) is the same as in
(\ref{new3-1}), the analysis in the proof of Theorem
\ref{microc-w} can be used to show the right hand side of
(\ref{new3-1-nn-f}) can be controlled by $\phi+ |\nabla
\phi|-C\sum_{i,j\in B}|\nabla W_{ij}|$. The theorem follows the
same argument as in the end of the proof of Theorem
\ref{microc-w-para}. \qed

\medskip

We now use Proposition \ref{thmw2-pro} to prove Theorem
\ref{thmw2-flow}. In fact, the local convexity condition on $F$ in
that theorem can be weakened to condition (\ref{eq2-4}).

\begin{theorem}\label{thmw2-flow-1}
Suppose $F(A,X, \vec n, t)$ is elliptic in $A$ and satisfies
condition (\ref{eq2-4}) for each fixed $\vec n\in \mathbb S^n$,
$t\in [0,T]$ for some $T>0$. Let $M(t)\subset \mathbb R^{n+1}$ be
compact hypersurface and it is a solution of (\ref{flow1}). If
$M_0$ is convex, then $M(t)$ is strictly convex for all $t\in
(0,T)$.
\end{theorem}

\noindent{\bf Proof of Theorem \ref{thmw2-flow-1}.} First, we may
approximate $M_0$ by a strictly convex $M^{\epsilon}_0$. By
continuity, there is $\delta>0$ (independent of $\epsilon$), such
that there is a solution $M^{\epsilon}(t)$ to (\ref{flow1}) with
$M^{\epsilon}(0)= M^{\epsilon}_0$ for $t\in [0,\delta]$. We argue
that $M^{\epsilon}(t)$ is strictly convex for $t\in [0,\delta]$.
If not, there is $t_0>0$, $M^{\epsilon}(t)$ is strictly convex for
$0\le t<t_0$, but there is one point $x_0$ such that
$(h_{ij}(x_0,t_0))$ is not of full rank. This is contradiction to
Proposition \ref{thmw2-pro}. Taking $\epsilon \to 0$, we conclude
that $M(t)$ is convex for all $t\in [0,\delta]$. This implies that
the set $t$ where $M(t)$ is convex is open. It is obviously
closed. Therefore, $M(t)$ is convex for all $t\in [0,T]$. Again,
by Proposition \ref{thmw2-pro}, $M(t)$ is strictly convex for all
$t\in (0,T]$. \qed

\begin{remark} If $n=2$, by Lemma \ref{homog2}, if $F(A)$ is homogeneous of
degree $k$ for either $k\ge 1$ or $k\le 0$, then $F$ satisfies
condition (\ref{eq2-4}) automatically.\end{remark}

\medskip

Let $(M,g)$ be a Riemannian manifold (not necessary compact), a
symmetric $2$-tensor $W$ is called a Codazzi tensor if $w_{ijk}$ is
symmetric with respect to indices $i,j,k$ in local orthonormal
frames. One of the important example of the Codazzi tensor is the
second fundamental form of hypersurfaces.

\begin{theorem}\label{codazzi}
Let $F(A,x)$ is elliptic and $F(A^{-1},x)$ is locally convex in
$(A,x)$. Suppose $(M,g)$ is a connected Riemannian manifold of
nonnegative sectional curvature, and $W$ is a semi-positive
definite Codazzi tensor on $M$ satisfying equation
\begin{equation}\label{ecod11}
F(g^{-1}W, x)=0 \quad \text{on $M$,}
\end{equation}
then $W$ is of constant rank and its null space is parallel.\end{theorem}

\noindent{\bf Proof.} Since the proof is
similar to the proof of Theorem \ref{microc}, we only indicate
some necessary modifications.

We use the same notations as in the proof of Theorem \ref{microc}.
As before, we set $\phi(x)
=\sigma_{l+1}(W(x))+\frac{\sigma_{l+2}(W(x))}{\sigma_{l+1}(W(x))}$
as in (\ref{def-phi}). As before, we want to establish corresponding
differential inequality (\ref{2-6}) in this case for the Codazzi
tensor $W$. We note that all the analysis in Section 3 carry through
without any change if we use local orthonormal frames, except the
commutators of derivatives. Since $W$ is Codazzi, we only need to
take care of commutators like $W_{\alpha \alpha,\beta
\beta}-W_{\beta \beta,\alpha \alpha}$. The Ricci identity states
\begin{eqnarray}\label{Ricci}
W_{\alpha \alpha, \beta\beta}=W_{\beta\beta, \alpha \alpha}
+R_{\alpha \beta \alpha \beta}
(W_{\alpha\alpha}-W_{\beta\beta}),
\end{eqnarray}
where $R_{\alpha \beta \alpha \beta}$ the sectional curvatures of
$(M,g)$. The assumption of nonnegativity of $R_{\alpha \beta
\alpha \beta}$ gives us a good sign, following the same lines of
the proof of Theorem \ref{microc-w}, we have the corresponding
differential inequality
\begin{equation}\label{2-6-1}
\sum_{\alpha \beta} F^{\alpha\beta}\phi_{\alpha\beta}(x)\leq
C_1(\phi(x)+|\nabla\phi(x)|)-\sigma_l(G)\sum_{\alpha \in G,
\beta\in B} F^{\alpha\alpha}R_{\alpha\beta
\alpha\beta}W_{\alpha\alpha}-C_2\sum_{i,j\in B}|\nabla W_{ij}|.
\end{equation}
The strong maximum principle implies $\phi\equiv 0$ in $M$, so $W$ is of
constant rank $l$. Again, by (\ref{2-6-1}), $\sum_{i,j\in B}|\nabla W_{ij}|\equiv 0$,
so the null space of $W$ is parallel. \qed
\medskip

\noindent{\bf Proof of Theorem \ref{codazzi-1}.} We deal with case
(2) of theorem first. Let $c=\min_{x\in M} W_s(x)$, where $W_s(x)$
is smallest eigenvalue of $W$ at $x$. We set $\tilde W=g^{-1}(W-c
g)$. Then $\tilde W$ is also a Codazzi tensor, it's rank is strictly
less than $n$ at some point, and it satisfies
\begin{equation} \tilde F(\tilde W)=F(g^{-1}\tilde W +cI)=constant.\end{equation}
By our assumption, $c\ge 0$, it follows from Corollary \ref{cor1c}
that $\tilde F$ satisfies condition (\ref{cond-c}). For $\phi(x)
=\sigma_{l+1}(\tilde W(x))+\frac{\sigma_{l+2}(\tilde
W)}{\sigma_{l+1}(\tilde W(x))}$, inequality (\ref{2-6-1}) is
valid. Therefore it follows from the proof of Theorem
\ref{microc-w}, $\phi\equiv 0$ in $M$. Now back to (\ref{2-6-1}),
the left hand side is identical to $0$, so is the right hand side.
By the assumption, $R_{\alpha\beta \alpha\beta}>0$ at some point.
It follows $G$ must be empty, that is $\tilde W\equiv 0$.

We now consider case (1), we follow the arguments in the proof of
Theorem \ref{n2} and Remark \ref{rkn2}. Let $\tilde W$ defined as before
($c$ may not necessary nonnegative in this case). $\tilde W$ is a
semi-positive Codazzi tensor, it's minimal rank $l$ is strictly
less than $2$ at some point, and it satisfies $\tilde F(\tilde
W)=F(g^{-1}\tilde W +cI)=0$, and $\tilde F$ is elliptic. If $l=0$,
the proof for case (2) carry through without change. If $l=1$,
i.e. $|G|=1$. At the given point, we may assume $\tilde W$ is
diagonal and $n\in G$. Differentiate equation $\tilde F(\tilde
W)=0$, as in the proof of Theorem \ref{n2}, we get
$$ \nabla
\tilde W_{nn}=O(\sum_{i,j\in B}\nabla \tilde W_{ij}).$$ Therefore,
$\nabla \tilde W_{nn}$ can be controlled. It follows from the
proof of Theorem \ref{microc-w}, inequality (\ref{2-6-1}) is
valid. In turn, we get $\phi\equiv 0$ in $M$. As in case (2),
since $R_{\alpha\beta \alpha\beta}>0$ at some point, we must have
$\tilde W\equiv 0$. \qed

\medskip

\begin{remark}
In spirit, our results are similar to
Hamilton's strong maximum principle \cite{Hamilton} for tensor equation
\begin{equation}\label{h1}  W_t=\Delta W+\Phi(W),\end{equation}
under the assumption that $V^T\Phi(W)V\ge 0$ for any null direction of $W$.
Our cases are different in the setting. For example, in the case of Theorem \ref{microc-w-para},
$W=(\nabla^2 u)$ satisfies
\begin{equation}\label{W1} W_t=F^{ij}\nabla_i\nabla_j W
+\Phi(\nabla W, W, \nabla u, u, x, t),\end{equation}
where $\Phi$ involves $\nabla W, W, \nabla u, u, x, t$.
Our main analysis is to
show $\Phi $ is controlled by $\phi+|\nabla \phi|$ near the null set of $\phi$.
\end{remark}

\begin{remark} Let $\lambda_{min}(t)=\min_{x\in M(t)} \{\mbox{smallest eigenvalue of h(x,t)}\}$.
If $F$ in (\ref{flow1}) is nonnegative and it depends only on $A$,
using Corollary \ref{cor1c} and (\ref{flow2}), by considering
$W=(h^i_j(x,t))-\lambda_{min}(s)I$, if $W$ has zero eigenvalue at
some time $t>s$, our argument in the above can show
\begin{equation}\label{2-6-1p}
\sum_{\alpha \beta}
F^{\alpha\beta}\phi_{\alpha\beta}(x)-\phi_t\leq
C_{1}\phi(x)+C_{2}|\nabla\phi(x)|-\sigma_l(G)\sum_{\alpha \in G,
\beta\in B} F^{\alpha\alpha}R_{\alpha\beta
\alpha\beta}W_{\alpha\alpha}.
\end{equation}
By Theorem \ref{thmw2-flow} the sectional curvature of $M(t)$ is
strictly positive, therefore the last term in (\ref{2-6-1p}) must
be vanishing, that is $W\equiv0$. In turn, Theorem
\ref{thmw2-flow} can be strengthened as follow:
\[\lambda_{min}(t)\ge \lambda_{min}(s), \quad \forall 0\le s\le
t\le T,\] if equality holds for some $s<t_0$, then
$(h^i_j(x,t))=\lambda_{min}(s)I$ is constant for all $s\le t$ and
for all $x$, that is $M(t)$ is a sphere for all $t\ge s$.
\end{remark}

\begin{remark}
Applying the same argument as in Remark \ref{wconds}, we can weaken local
convexity condition on $F$ in
Theorem \ref{codazzi-1} and Theorem \ref{codazzi}. Let
\[\mathcal{D}_{W(x)}=\{\mbox{$r$ diagonal} |\quad r=Qg^{-1}(x)W(x)Q^T   \mbox{for some $Q\in
O(n)$}\},\]
\[\tilde {D}_{W(x)}^{\delta}=\{A| \quad \|A^{-1}-r\|\le \delta,
\mbox{ for some $r\in \mathcal{D}_{u(x)}$}\}.\]
In this case, we only need the condition: there is $\delta>0$,
\begin{equation}\label{newccc}
\mbox{ $F(A^{-1},x)$ is locally convex in $\tilde
{D}_{W(x)}^{\delta}\times \mathcal O$ .}\end{equation} We note
that when $M$ is compact, for given Codazzi tensor $W$ on $M$,
there is $\lambda>0$ such that $\tilde W=\lambda g-W\ge 0$
everywhere. If $F(W)$ is concave in $W$, then $\tilde
F(g^{-1}\tilde W)=-F(\lambda I-g^{-1}\tilde W)$ satisfies
condition (\ref{newccc}).
\end{remark}


\begin{thebibliography}{99}

\bibitem{Alex1} A.D. Alexandrov, {\em Zur Theorie der gemischten Volumina von
konvexen korpern, III. Die Erweiterung zweeier Lehrsatze
Minkowskis uber die konvexen polyeder auf beliebige konvexe
Flachen ( in Russian)} Mat. Sbornik N.S. \textbf{3}, (1938),
27-46.

\bibitem{Alex} A. V. Alexandrov, {\em Über konvexe Flächen mit ebenen Schattengrenzen},
(Russian) Rec. Math. N. S. [Mat. Sbornik] {\bf 5(47)}, (1939),
309--316.

\bibitem{Alex2}
A.D. Alexandrov, {\em Uniqueness theorems for surfaces in the
large. I} (Russian), Vestnik Leningrad. Univ. \textbf{11} (1956),
5--17. English translation: AMS Translations, series 2,
\textbf{21}, (1962), 341-354.

\bibitem{All97}
O. Alvarez, J.M. Lasry and P.-L. Lions, {\em Convexity viscosity
solutions and state constraints}, J. Math. Pures Appl. {\bf 76},
(1997), 265-288.
\bibitem{Andrews} B. Andrews, {\em Pinching estimates and motion of
hypersurfaces by curvature functions}. J. Reine Angew. Math. {\bf
608}, (2007), 17--33.

\bibitem {BL}
H.J. Brascamp and E.H. Lieb, {\em On extensions of the
Bruun-Minkowski and Prekopa-Leindler theorems, including
inequalities for log-concave functions, with an application to the
diffusion equation}, J. Funct. Anal., {\bf 22}, (1976), 366-389.

\bibitem {CF85}
L. Caffarelli and A. Friedman, {\em Convexity of solutions of some
semilinear elliptic equations}, Duke Math. J. {\bf 52},  (1985),
431-455.

\bibitem{CGM} L. Caffarelli, P. Guan and X. Ma, {\em A constant rank theorem for
solutions of fully nonlinear elliptic equations}, Communications
on Pure and Applied Mathematics, \textbf{60}, (2007), 1769-1791 .

\bibitem{CNS}
L. Caffarelli, L. Nirenberg and J. Spruck, {\em The Dirichlet
problem for nonlinear second order elliptic equations,  III:
Functions of the eigenvalues of the Hessian}, Acta Math.,
\textbf{155}, (1985), 261-301.

\bibitem{CS} L. Caffarelli and J. Spruck, {\em Convexity
properties of solutions to some classical variational problems},
Comm. in Partial Differential Equations, {\bf 7}, (1982),
1337-1379.

\bibitem{CY}
S. Y. Cheng and S. T. Yau, {\em Hypersurfaces with constant scalar
curvature}, Math. Ann. {\bf 225} (1977), 195--204.

\bibitem{Chern} S. S. Chern, {\em Some new characterizations of the Euclidean
sphere}, Duke Math. J. {\bf 12}, (1945). 279--290.

\bibitem{E-H} K. Ecker and G. Huisken, {\em Immersed hypersurfaces with
constant Weingarten curvature}, Math.Ann. {\bf 283}(1989), 329-332.

\bibitem{Guan-Duke} P. Guan, {\em $C^2$ A Priori Estimates for Degenerate
Monge-Ampere Equations}, Duke Mathematical Journal, \textbf{86},
(1997), 323-346.

\bibitem{GLZ} P. Guan, Q. Li and X. Zhang, {\em A uniqueness theorem in K\"ahler geometry},
preprint, 2007.

\bibitem{GLM1} P. Guan, C.S. Lin and X.N. Ma, {\em The
Christoffel-Minkowski problem II: Weingarten curvature equations},
Chinese Annals of Mathematics, Series B.\textbf{27}, (2006), 595-614.

\bibitem{gm}
P. Guan and X.N. Ma, {\em The Christoffel-Minkowski Problem I:
Convexity of Solutions of a Hessian Equations}, Inventiones Math.,
{\bf 151}, (2003), 553-577.

\bibitem{GMZ} P. Guan, X.N. Ma and F. Zhou,
{\em The Christoffel-Minkowski problem III: existence and
convexity of admissible solutions},  Comm. Pure and Appl. Math.
{\bf 59}, (2006), 1352-1376.

\bibitem{Hamilton} R.S. Hamilton, {\em Four manifolds with positive curvature operator},
J. Differential Geometry,
{\bf 24}, (1986), 153-179.

\bibitem{HN59} P. Hartman and L. Nirenberg, {\em On spherical image maps
whose Jacobians do not change sign}. Amer. J. Math. \textbf{81},
(1959), 901--920.

\bibitem{Huisken} G. Huisken, {\em Flow by mean curvature of convex surfaces into spheres}.
J. Differential Geometry, {\bf 20} (1984), 237-266.

\bibitem{HS99} G. Huisken and C. Sinestrari, {\em Convexity estimates for mean
curvature flow and singularities of mean convex surfaces}. Acta
Math. {\bf 183}, (1999), 45-70.

\bibitem{Ka86}
B. Kawohl, {\em A remark on N.Korevaar's concavity maximum
principle and on the asymptotic uniqueness of solutions to the
plasma problem}, Math. Methods Appl. Sci., {\bf 8}, (1986),
93-101.


\bibitem{Ke85}
A.U. Kennington, {\em Power concavity and boundary value
problems}, Indiana Univ. Math. J., {\bf34}, (1985),  687-704.

\bibitem{Ko831}
N.J. Korevaar, {\em Capillary surface convexity above convex
domains}, Indiana Univ. math. J., {\bf 32}, (1983), 73-81.

\bibitem{Ko832}
N.J. Korevaar, {\em Convex solutions to nonlinear elliptic and
parabolic boundary value problems}, Indiana Univ. math. J., {\bf
32}, (1983),  603-614.


\bibitem{Kl87}
N.J. Korevaar and J. Lewis, {\em Convex solutions of certain
elliptic equations have constant rank hessians}, Arch. Rational
Mech. Anal. {\bf 91}, (1987),  19-32.

\bibitem{SWYY} I. Singer, B. Wong, S.T. Yau and Stephen S.T. Yau, {\em
An estimate of gap of the first two eigenvalues in the Schrodinger
operator }, Ann. Scuola Norm. Sup. Pisa Cl. Sci.(4), {\bf 12}
(1985), 319-333.

\bibitem{T} F. Treves, {\em A new proof of the subelliptic estimates},
Comm. Pure Appl. Math. {\bf 24}, (1971), 71-115.




\end{thebibliography}
\end{document}